\documentclass{gtart}

\def\ifplaintex{\expandafter\ifx\csname documentclass\endcsname\relax}

\def\gtp{{\mathsurround=0pt\it $\cal G\mskip-2mu$eometry \&\ 
$\cal T\!\!$opology $\cal P\!$ublications}}  % GT publications

\def\recd{{\small Received:\qua\receiveddate\ifx\reviseddate\relax
\else\qquad Revised:\qua\reviseddate\fi\par}} 

\def\lognumber#1{\def\thelognumber{#1}}
\def\volumenumber#1{\def\thevolumenumber{#1}}
\def\volumeyear#1{\def\thevolumeyear{#1}}
\def\papernumber#1{\def\thepapernumber{#1}}
\def\pagenumbers#1#2{\def\startpage{#1}\def\finishpage{#2}}
\def\published#1{\def\publishdate{#1}}

\def\received#1{\def\receiveddate{#1}}
\def\revised#1{\def\reviseddate{#1}}
\def\accepted#1{\def\accepteddate{#1}}

\def\asciiaddress#1{\def\theasciiaddress{#1}}

\long\def\asciiabstract#1{\long\def\theasciiabstract{#1}}

\let\\\par\let\thelognumber\relax\let\thevolumenumber\relax
\let\thepapernumber\relax\let\thevolumeyear\relax\let\startpage\relax
\let\finishpage\relax\let\publishdate\relax\let\receiveddate\relax
\let\reviseddate\relax\let\accepteddate\relax\let\theasciititle\relax
\let\theasciiauthors\relax\let\theasciiaddress\relax
\let\theasciiabstract\relax

\let\theasciiemail\relax

\ifplaintex
\font\logobig=cmssbx10 scaled 3836
\font\logomed=cmssbx10 scaled 2557
\else
\font\logobig=cmssbx10 scaled 4200
\font\logomed=cmssbx10 scaled 2800
\fi

\long\def\makeagttitle{   %%% start of definition of \makeagttitle
\count0=\startpage
\agt\hfill      %   Journal title (top left) 
\hbox to 45truept{\vbox to 0pt{\vglue -13truept{\logomed A\kern -.37em{\logobig 
T}\kern -.38em G}\vss}\hss}
\break
{\small Volume \thevolumenumber\ (\thevolumeyear)
\startpage--\finishpage\nl
Published: \publishdate}

\vglue .25truein

{\parskip=0pt\leftskip 0pt plus
1fil\def\\{\par\smallskip}{\Large\bf\thetitle}\par\medskip} \vglue
0.05truein

{\parskip=0pt\leftskip 0pt plus 1fil\def\\{\par}{\sc\theauthors}
\par\medskip}%
 
\vglue 0.03truein 

{\small\leftskip 25truept\rightskip 25truept{\bf Abstract}\stdspace\theabstract

{\bf AMS Classification}\stdspace\theprimaryclass
\ifx\thesecondaryclass\relax\else; \thesecondaryclass\fi\par
{\bf Keywords}\stdspace \thekeywords\par}\vglue 7truept

}   %%%% end of definition of \makeagttitle

\ifplaintex
\font\phead=cmsl9 scaled 950
\font\pnum=cmbx10 scaled 913
\font\pfoot=cmsl9 scaled 950
\headline{\vbox to 0pt{\vskip -4.5mm\line{\small\phead\ifnum
\count0=\startpage ISSN 1472-2739 (on-line) 1472-2747 (printed)
\hfill {\pnum\folio}\else\ifodd\count0\def\\{ }% 
\ifx\theshorttitle\relax\thetitle\else\theshorttitle\fi\hfill{\pnum\folio}
\else\def\\{ and }{\pnum\folio}\hfill\ifx\theshortauthors\relax\theauthors
\else\theshortauthors\fi\fi\fi}\vss}}
\footline{\vbox to 0pt{\vglue 0mm\line{\small\pfoot\ifnum\count0=\startpage
\copyright\ \gtp\hfill\else
\agt, Volume \thevolumenumber\ (\thevolumeyear)\hfill\fi}\vss}}
\else
\font=cmsl9 scaled 1050
\font\lnum=cmbx10 
\font=cmsl9 scaled 1050
\makeatletter
\def\@oddhead{{\small\lhead\ifnum\count0=\startpage ISSN 1472-2739 
(on-line) 1472-2747 (printed)\hfill {\lnum\number\count0}\else\ifodd\count0
\def\\{ }\ifx\theshorttitle\relax \thetitle \else\theshorttitle\fi\hfill
{\lnum\number\count0}\else\def\\{ and }{\lnum\number\count0}
\hfill\ifx\theshortauthors\relax 
\theauthors\else\theshortauthors\fi\fi\fi}}\def\@evenhead{\@oddhead}
\def\@oddfoot{\small\lfoot\ifnum\count0=\startpage\copyright\ \gtp\hfill\else
\agt, Volume \thevolumenumber\ (\thevolumeyear)\hfill\fi}
\def\@evenfoot{\@oddfoot}
\makeatother
\fi
\let\maketitlepage\makeagttitle
\let\makeshorttitle\maketitlepage
\let\maketitle\maketitlepage

\newwrite\gtoutfile
\long\gdef\makeheadfile{  %%% start of definition of \makeheadfile
{\def\\{, }\def\s{ }
\immediate\openout\gtoutfile head.xxx
\immediate\write\gtoutfile{To: math@arxiv.org}
\immediate\write\gtoutfile{Subject: put OR rep NNNNN:ppppp}
\immediate\write\gtoutfile{--text follows this line--}
\immediate\write\gtoutfile{Proxy-for: \ifx\theasciiauthors\relax
\theauthors\else\theasciiauthors\fi\s<\ifx\theasciiemail\relax\theemail\else\theasciiemail\fi>}
\immediate\write\gtoutfile{\noexpand\\}
\immediate\write\gtoutfile{Authors: \ifx\theasciiauthors\relax
\theauthors\else\theasciiauthors\fi}
{\def\\{ }\immediate\write\gtoutfile{Title: \ifx\theasciititle\relax
\thetitle\else\theasciititle\fi}}
\immediate\write\gtoutfile{Subj-class: GT or SG, GR etc}
\immediate\write\gtoutfile{MSC-class: \theprimaryclass\ifx\thesecondaryclass\relax\else, \thesecondaryclass\fi}
\immediate\write\gtoutfile{Journal-ref: Algebr. Geom. Topol. \thevolumenumber\s
(\thevolumeyear) \startpage-\finishpage}
\immediate\write\gtoutfile{Comments: Published by Algebraic and
Geometric Topology at}
\immediate\write\gtoutfile{\s\s\s  http://www.maths.warwick.ac.uk/agt/AGTVol\thevolumenumber/agt-\thevolumenumber-\thepapernumber.abs.html}
\immediate\write\gtoutfile{\noexpand\\}
\immediate\write\gtoutfile{}
\ifx\theasciiabstract\relax
\immediate\write\gtoutfile{\theabstract}\else
\immediate\write\gtoutfile{\theasciiabstract}\fi
\immediate\write\gtoutfile{}
\immediate\write\gtoutfile{\noexpand\\}
\immediate\write\gtoutfile{}
\immediate\closeout\gtoutfile}}  %%% end of definition of \makeheadfile

\def\maketitlepage{\makeagttitle\makeheadfile}
\let\makeshorttitle\maketitlepage
\let\maketitle\maketitlepage

\lognumber{15}
\volumenumber{3}
\volumeyear{2003}
\papernumber{15}
\published{5 June 2003}
\pagenumbers{473}{509}
\received{26 September 2002}
\revised{12 May 2003}
\accepted{5 June 2003}

\usepackage{amsmath,diagrams,amssymb}

\newtheorem{Thm}{Theorem}[section]
\newtheorem{Lem}[Thm]{Lemma}
\newtheorem{Cor}[Thm]{Corollary}
\newtheorem{Prop}[Thm]{Proposition}
\theoremstyle{definition}
\newtheorem{Def}[Thm]{Definition}
\newtheorem{Rem}[Thm]{Remark}

\begin{document}

\title{Transfer and complex oriented cohomology rings}

\author{Malkhaz Bakuradze\\Stewart Priddy}
\address{Razmadze Mathematical Institute, Tbilisi 380093, Republic of Georgia 
\\{\rm and}\\Department of Mathematics, Northwestern University, Evanston, 
IL 60208, USA}
\asciiaddress{Razmadze Mathematical Institute, Tbilisi 380093, Republic of Georgia 
\\and\\Department of Mathematics, Northwestern University, Evanston, 
IL 60208, USA}

\email{maxo@rmi.acnet.ge, priddy@math.northwestern.edu}

\begin{abstract}
For finite coverings we elucidate the interaction between transferred
Chern classes and Chern classes of transferred bundles.  This involves
computing the ring structure for the complex oriented cohomology of
various homotopy orbit spaces.  In turn these results provide
universal examples for computing the stable Euler classes (i.e.\
$Tr^*(1)$) and transferred Chern classes for $p$-fold covers.
Applications to the classifying spaces of $p$-groups are given.
\end{abstract}
\asciiabstract{%
For finite coverings we elucidate the interaction between transferred
Chern classes and Chern classes of transferred bundles.  This involves
computing the ring structure for the complex oriented cohomology of
various homotopy orbit spaces.  In turn these results provide
universal examples for computing the stable Euler classes
(i.e. Tr^*(1)) and transferred Chern classes for p-fold
covers. Applications to the classifying spaces of p-groups are
given.}



\maketitle

\section{ Introduction}

For various examples of finite groups the complex oriented cohomology ring
coincides with its subring generated by Chern classes
\cite{SY},\cite{TY}, \cite{Y}. Even more groups are good in the sense that
their Morava $K$-theory is generated by transferred Chern classes
of complex representations of subgroups \cite{HKR}. Special effort was needed
to find an example of a group not good in this sense \cite{Kr}. Thus the
relations in the complex oriented cohomology ring of a finite group derived
from formal properties of the transfer should play a major role. The purpose
of  this paper is to elucidate for finite coverings the interaction between
transferred Chern classes and Chern classes of transferred bundles.

Let $p$ be a prime and let $G \leq \Sigma_p$ be a subgroup of the symmetric
group.
In this paper we consider the complex oriented cohomology
of homotopy orbit spaces $X^p_{hG} = EG {\times}_G X^p$.
Several authors have computed these cohomology groups,
\cite{MS}, \cite{Hun}, \cite{Hun1}, \cite{HKR}, however
we are particularly interested in the
ring structure and thereby explicit formulas for the transfer.
Thus we are led to consider Fibrins reciprocity, the relation between
cup products and transfer:
$$
Tr^*(x)y = Tr^*(x\rho^*(y))
$$
(formula (i) of Section 2) where $\rho: EG {\times} X^p \to X^p_{hG}$ is the covering
projection
and
$$
  Tr^*: E^*(X^p) \to E^*(X^p_{hG})
$$
is the associated transfer homomorphism.

Let $\pi \leq \Sigma_p$ be the subgroup of cyclic permutations of order $p$.
Our results for $MU^*(X^p_{h\pi}),$ $X = \mathbf{C}P^{\infty}$, and
$\xi \rightarrow \mathbf{C}P^{\infty}$ the canonical complex line bundle,
provide a universal example which enables us to write explicitly the  Chern
classes $c_1,\ldots,c_{p-1}$ of the transferred bundle $\xi_{\pi}$
as certain formal power series in the Euler class $c_p(\xi_{\pi})$ with
coefficients in $E^*(B\pi)$ plus certain transferred classes of the
bundle $\xi$. In Section 3 we give an algorithm for computing these coefficients.

In paricular for $E=BP$, Brown-Peterson cohomology, the
coefficients of this formal power series are invariant under the action
of the normalizer of $\pi$ in $\Sigma_p$. This enables us to give the similar
results for $\Sigma_p $ coverings. Moreover in Section 4 we compute the algebra
$BP^*(X^p_{h\Sigma_p})$ and show that its multiplicative
structure is completely determined by Frobenius reciprocity.

In addition for $E=K(s)$, Morava $K$-theory, the computations
become easier: we show in Section 5 that the formal power series in the
algorithm above descend to polynomials. We derive an alternative way for
calculation and give some examples.

Section 6 is devoted to extending some results of \cite{HKR} in Morava
$K$-theory. In particular we show that if $X$ is good then $X^p_{h\pi}$ is good.

In Section 7 we give some applications to classifying spaces of finite groups.

We would like to thank several people for their help during the course of this
work: D. Ravenel for supplying us with the proof of Lemma \ref{LemF1}, M. Jibladze for some
Maple programs used in the examples of Section 5, and finally the referee who
suggested many improvements.

The first author was supported by CRDF grant
GM1 2083 and by the Max-Planck-Institut f\"ur Mathematik.

A word about notation: in Sections 2, 3, 4 and 7 we denote
$\mathbf{C}P^{\infty}$ simply by $X$.

\section{ Preliminaries}

We recall that a multiplicative cohomology theory $E^*$ is called {\it complex oriented} if
there exists a Thom class, that is, a class $u \in E^2(\mathbf{C}P^{\infty})$ that restricts
to a generator of the free one-dimensional $E^*$ module $E^2(\mathbf{C}P^1)$. The universal
example is complex cobordism $MU^*$. Then
$$
E^*(\mathbf{C}P^{ \infty})=E^*[[x]],
$$
where $x$ is the Euler class of the canonical complex line
bundle $ \xi$ over $\mathbf{C}P^{ \infty} = BU(1)$. Further
$$
E^*(BU(1)^p)=E^*[[x_1,\ldots,x_p]],
$$
where $x_i = c_1( \xi_i)$ and $ \xi_i$ is the pullback bundle over $BU(1)^p$
by the projection $BU(1)^p \rightarrow BU(1)$ on the $i$-th factor.

Much of our paper is written in terms of transfer maps \cite{Ada, KP}
and formal group laws. Let us give a brief review of formal properties of the
transfer. For a finite covering
$$
\rho:X \rightarrow X/G
$$
there is a stable transfer map
$$
  Tr = Tr(\rho):X/G^+ \rightarrow X^+ .
$$
For any multiplicative cohomology theory $E^*$, {\it Frobenius reciprocity} holds i.e.,
the induced map $Tr^*$ is a map of
$E^*(X/G)$ modules

(i)\qua $Tr^*(x \rho^*(y))=Tr^*(x) y, \: x \in E^*(X), \: y \in E^*(X/G).$

\noindent For example

(ii)\qua $Tr^*( \rho^*(y))=Tr^*(1) y.$

The element $Tr^*(1) \in E^0 (X/G)$ is called the {\it index or
stable Euler class } of the covering $ \rho$. The following additional properties
of the transfer will be used:

(iii)\qua The transfer is natural with respect to pullbacks;

(iv)\qua $Tr( \rho_1 \times \rho_2) = Tr( \rho_1) \wedge Tr( \rho_2) $;

(v)\qua If $ \rho = \rho_2 \rho_1 $, then $Tr( \rho)=Tr( \rho_2)Tr( \rho_1)$.

More generally for a covering projection
$$
   \rho_{H,G}: X/H \rightarrow X/G
$$
with $H \leq G$ there is a stable transfer map
$$
  Tr_{H,G}: X/G^+ \rightarrow X/H^+.
$$
To ease
notation if  $H = e$, as above, we write projection and transfer in
equivalent ways $\rho = \rho_G$,
$Tr = Tr(\rho) = Tr_G$.

The reverse composition to (ii) is given by:

(vi)\qua (Double coset formula)\qua If $K,H \leq G$ then
$$
\rho_{K,G}^* Tr_{H,G}^* = \sum_x Tr_{K\cap H^{x},K}^*\circ {x^{-1}}^*\circ \rho_{K^{x^{-1}}\cap H ,H}^*
$$
where the sum is taken over a set of double coset representatives $x \in K\backslash G/H$.
Here $H^x = xHx^{-1}$.

For a regular covering $ \rho_{H,G}$, i.e.\ $H \trianglelefteq G$,
$$
\rho_{H,G}^* Tr_{H,G}^*(x)= N(x) = \sum_{g \in G/H} g^* (x),
$$
where $N(x)$ is called the {\it norm} or {\it trace} of $x$.

In subsequent sections the reduced transfer $Tr_{H,G}: X/G \rightarrow X/H$
is used.

We recall  Quillen's formula \cite {Q, Die}.
First,
$$
E^*(B\mathbf{Z}/p)=E^*[[z]]/([p](z)),
$$
where  $x$ is the Euler class of a faithful one-dimensional complex representation
of $\mathbf{Z}/p$ and $[p](z)$ is the $p$-series or $p$-fold iterated formal sum.
Then
\begin{equation}
Tr_{\mathbf{Z}/p}^*(1)=[p](z)/z, \label{eq: qq}
\end{equation}
where  $Tr_{\mathbf{Z}/p}^*$ is the transfer homomorphism for the universal $\mathbf{Z}/p$-covering
$E\mathbf{Z}/p \rightarrow B\mathbf{Z}/p$.
The relation $[p](z)=0$ is equivalent to the transfer
relation
$$
zTr_{\mathbf{Z}/p}^*(1)=Tr_{\mathbf{Z}/p}^*(c_1(\mathbf{C}))=Tr_{\mathbf{Z}/p}^*(0)=0
$$
obtained by applying (ii).
Of course since the transfer is natural, Quillen's formula enables us
to compute the stable Euler class for any regular $\mathbf{Z}/p $ covering.

In this spirit, let
$$
\pi = \langle t \rangle \leq \Sigma_p
$$
be the subgroup of cyclic permutations of order $p$. For a given free action
of $\pi$  on a space $Y$ with a given complex line
bundle $ \eta \rightarrow Y $ we have an equivariant map
$$
\eta_{ \pi}=(g_1,\ldots,g_p):Y \rightarrow BU(1)^p,
$$
where $g_i$ classifies the line bundle $t^{i-1} \eta$.

So by naturality of the transfer, the computation of transferred Chern classes
$Tr^*(c_1^i( \eta))$, $i \geq 1$ for cyclic coverings can be reduced to
the covering
$$
 \rho_{ \pi}: E\pi \times (BU(1))^p  \rightarrow E\pi \times_{\pi}(BU(1))^p,
$$
as the universal example.

Similarly for the symmetric group.

Let $ \xi $ be the canonical complex line bundle
over $\mathbf{C}P^{ \infty} =BU(1)$ and $ \xi_i$ be the pullback bundle over $BU(1)^p$
by the projection on the $i$-th factor as before. Then
$MU^*(BU(1)^p)=MU^*[[x_1,\ldots,x_p]]$, $x_i = c_1( \xi_i)$ and $x_1 \cdots x_p$
is the Euler class of the bundle $ \xi^{ \times p}= \oplus \xi_i$.

Note that by transfer property (v), $Tr( \rho_{ \pi})^* $ has the same value on
the Chern classes $x_1,\ldots,x_p$:
the group $ \pi$ permutes the $x_i$ and $\rho_{ \pi}t= \rho_{ \pi}$,
$t \in \pi$. Thus in computations of the transfer we sometimes write these Chern classes
in an equivalent way $x, tx, \ldots, t^{p-1}x$.

For the sphere bundle $S( \xi^{ \times p}),$
one has
\begin{equation}
MU^*(S(\xi^{ \times p}))=MU^*[[x_1,\ldots,x_p]]/(x_1 \cdots x_p).
\label{eq:133}
\end{equation}
Then for the trace map
$$N=1+t+\cdots +t^{p-1}$$
we have $ker N = Im(1-t)$, $t \in \pi$ in
$MU^*BU(1)^p$ and after restricting $N$ to $MU^*(S( \xi^{ \times p}))$
we have the exact sequence
\begin{equation}
 \cdots\!  \leftarrow MU^*(S( \xi^{ \times p})) \stackrel{N}{ \leftarrow}
MU^*(S( \xi^{ \times p}))
\stackrel{1-t}{ \leftarrow} MU^*(S( \xi^{ \times p})) \stackrel{N}{ \leftarrow}
MU^*(S( \xi^{ \times p}))\leftarrow\! \cdots 
\end{equation}
Then let $\xi_{ \pi} =E \pi \times_{ \pi} \xi^{ \times p}$ be the Atiyah
transfer bundle \cite{At},
\begin{equation}
S( \xi_{ \pi})=E\pi \times_{ \pi} S( \xi^{ \times p}) \label{eq:14}
\end{equation}
be its sphere bundle and
\begin{equation}
D( \xi_{ \pi})=E\pi \times_{ \pi} D( \xi^{ \times p})
\end{equation}
be its disk bundle. Let $X=\mathbf{C}P^{\infty}$ then $D( \xi_{ \pi})$ is homotopy equivalent to
$X^p_{h \pi} = B( \pi \wr U(1))$.

The cofibration $D( \xi_{ \pi })/ S( \xi_{ \pi})= (X^p_{h \pi})^{ \xi_{ \pi}} $
gives a long exact sequence
\begin{equation}
\cdots  \leftarrow MU^* (S( \xi_{ \pi})) \leftarrow MU^*(X^p_{h \pi})
\stackrel{ \times c_p}{ \leftarrow} MU^*((X^p_{h \pi})^{ \xi_{ \pi}}) \leftarrow \cdots  \label{eq:16}
\end{equation}
where $(X^p_{h \pi})^{ \xi_{ \pi}}$ is the Thom space of the bundle
$ \xi_{ \pi}$ and the right homomorphism is multiplication by the Euler class
$c_p =c_p ( \xi_{ \pi})$.

Since the diagonal of $ BU(1)^p $ is fixed under the
permutation action of $ \pi$, the inclusion
$E \pi \rightarrow E \pi \times BU(1)^p ; \: x \rightarrow (x,fixpoint)$
defines the inclusions
$$
i :B \pi \rightarrow X^p_{h \pi}
$$
$$
i_0 : B \pi \rightarrow S( \xi_{ \pi}).
$$
The projection
$ \varphi : X^p_{h\pi} \rightarrow B\pi$ induced by $ \pi \wr U(1) \to \pi$
defines the projection
\begin{equation}
\varphi_0 : S( \xi_{ \pi}) \rightarrow B\pi  \label{eq:17}
\end{equation}
and the compositions $ \varphi_0 i_0 $, $ \varphi i$ are the identity. We can
consider $S( \xi_{ \pi})$ as a bundle over $B\pi$ with fiber $S(\xi^{\times p})$.

Let $ \eta $ be the canonical line bundle over $B\pi$ and
$$
\theta = \varphi^*( \eta) \rightarrow X^p_{h \pi}
$$
be the pullback bundle. Thus  $i^*( \theta)= \eta $ and
$i^*( \xi_{ \pi})= \mathbf{C} + \eta +\cdots + \eta^{p-1}.$

Consider the pullback diagram
\[
\begin{diagram}
E\pi \times S( \xi^{ \times p})  & \rTo{} & E\pi \times BU(1)^p  \\
\dTo{\rho_0 } &  &  \dTo{ \rho }  \\
S( \xi_{ \pi})& \rTo{} & X^{p}_{h\pi}
\end{diagram}
\]
Let $Tr=Tr({ \rho}) $ be the transfer of the covering $ \rho $,
and
$Tr_0 : S( \xi_{ \pi}) \rightarrow S( \xi^{ \times p} ) $ the
transfer map of $ \rho_0$.

We will often refer to the following lemma which follows from (3)
and Frobenius reciprocity.

\begin{Lem} In $MU^*(X^p_{h\pi}),$ $ImTr^* \bigcap Ker(\rho^*)=0. $
\end{Lem}
\proof $\rho^*(Tr^*(a))=N(a)=0 \Rightarrow a \in Im(1-t) \Rightarrow Tr^*(a)=0$.\endproof

\begin{Rem} Lemma 2.1 is valid only in complex oriented cohomology $E^*$
with torsion free coefficient ring. This lemma is used in the proof
of Theorem 3.1 in complex cobordism and in the second statement of
Theorem 4.6 in Brown-Peterson cohomology. By naturality, these results
hold for all $E^*$ in the first case and all $p$-local $E^*$ in the second.
\end{Rem}

\section{Transferred Chern classes for cyclic coverings}

In this section we prove our main result for cyclic coverings, Theorem \ref{ThmB}.

In the notation of the previous section the $k$-th Chern class of the bundle
$\xi^{ \times p} $ is the elementary symmetric function $ \sigma_k(x_1,\ldots,x_p)$
in Chern classes $x_i$ and is the sum of $\binom pk$
elementary monomials. The action of $\pi$ on the set of these monomials gives
us $p^{-1}\binom pk$ orbits and the transfer homomorphism is constant on orbits
by transfer property (v) or (iii).

Let $E^*$ be a complex oriented cohomology theory. For $k = 1,\ldots,p-1$, let
$$\omega_k=\omega_k(x_1,\ldots,x_p) \in E^*(BU(1)^p)$$
be the sum of representative monomials one from each of these orbits. The value of
$Tr^* (\omega_k)$ does not depend on the choice of $\omega_k$ since $\omega_k$ is defined
modulo $Im(1-t)$ and on the elements of $Im(1-t)$ the transfer homomorphism
is zero again by (v). In other words we can take any $\omega_k$ for which
$N\omega_k= \sigma_k (x_1,\ldots,x_p)$ holds. As we shall explain in Corollary \ref{CorX}
of Theorem \ref{ThmJ}, the following result enables us to calculate the transfer
on all elements whose norm is symmetric.

For ease of notation let $X = \mathbf{C}P^{\infty}$ and $c_j = c_j(\xi_{\pi})$, 
$j = 1,\ldots,p$.

\begin{Thm}
\label{ThmJ}
We can construct explicit elements
$$
\delta_i^{(k)} \in \tilde{E}^*(B\pi), \;  k=1,\ldots,p-1,
$$
such that
$$
Tr^*(\omega_k)=c_k+ \sum_{i \geq0} \varphi^*( \delta_i^{(k)})c_p^i
$$
for the transfer of the covering $ \rho : X^p \rightarrow X^p_{h\pi}$.

\end{Thm}

Before constructing the elements $\delta_i^{(k)}$ in Section 3.2 we first prove their existence.

\subsection{ Complex cobordism of $(\mathbf{C}P^{\infty})^p_{h \pi}$ }

\begin{Thm}
\label{ThmB}
In $MU^*(X^p_{h\pi})$

{\rm(a)}\qua The annihilator of the Chern class $ c=c_1 ( \theta )$ coincides with
$ImTr^*$;

{\rm(b)}\qua Multiplication by  $c_p = c_p( \xi_{ \pi})$ is a monomorphism;

{\rm(c)}\qua Any element of $Ker(\rho^*)$ has the form
$ \sum_{k \geq0} \varphi^*( \delta_k) c_p^k $,
for some elements $ \delta_k \in  \tilde{MU}^*(B \pi)$.

{\rm(d)}\qua  For $\pi = \mathbf{Z}/2$,
$$
MU^*B(\pi \wr U(1))=MU^*[[c,c_1,c_2]]/(c_1-c^*_1, c_2-c^*_2)
$$
$$
=MU^*(B\pi)[[Tr^*(x), c_2]]/(cTr^*(x)),
$$
where $c_i=c_i( \xi_{ \pi})$, $c_i^*=c_i ( \xi_{ \pi} \otimes_{\mathbf{C}} \theta)$,
and $x$ are Chern characteristic classes with
$x\in MU^*(BU(1)^2) = MU^*[[x,tx]]$.
\end{Thm}

We need the following lemma.

\begin{Lem}
\label{LemH}
The left homomorphism in the long exact sequence (\ref{eq:16}) is an epimorphism
and thus gives a short exact sequence
$$
0 \leftarrow MU^* (S( \xi_{ \pi})) \leftarrow MU^*(X^p_{h\pi}) \stackrel{ \times c_p}{ \leftarrow}
MU^*(X^p_{h\pi})^{ \xi_{ \pi}} \leftarrow 0.
$$
Moreover there is a space $X_{\pi}$ and a stable equivalence
$$
\varphi_0 \vee f_{ \pi} : S( \xi_{ \pi}) \to B\pi \vee X_{ \pi} ,
$$
with $f_{ \pi} $ factoring through the following composite map
$$
S( \xi_{ \pi}) \rightarrow X^p_{h\pi} \stackrel{Tr}{ \rightarrow}
E \pi \times BU(1)^p
$$
and $ \varphi_0$ as in (\ref{eq:17}).
\end{Lem}

\proof Consider the Serre spectral sequence for the fibration
(\ref{eq:17})
$$
S( \xi^{ \times p}) \rightarrow S( \xi_{ \pi})
\stackrel { \varphi_0 } { \rightarrow} B \pi.
$$
$ E_2^{i,j}=H^i(\pi,H^j(S( \xi^{ \times p}); \mathbf{F}_q)) $ with the action of
$ \pi$ on
$H^*(S( \xi^{ \times p}); \mathbf{F}_q)$ by permutations of the cohomological Chern
classes.

When $q=p$,
$E_2^{0,j}= H^j (S( \xi^{ \times p}); \mathbf{F}_p)^{ \pi}$ and
$E_2^{i,0}=H^i (B\pi; \mathbf{F}_p).$

Then in positive dimensions
$H^*(S( \xi^{ \times p}); \mathbf{F}_q)= \mathbf{F}_q[x_1,\ldots,x_p]/(x_1 \cdots x_p) $
is a permutation representation of $ \pi$ acting on
monomials which have degree zero in at least one indeterminate.
This is a free $\mathbf{F}_q[\pi]$-module since all
the monomials that are fixed under this action have been factored out after
quotienting by the ideal $(x_1 \cdots x_p)$.
Hence the cohomology of $ \pi$ with coefficients in this module is
trivial in positive dimensions, i.e.\
$E_2^{i,j}=0$ when $i,j > 0$. Thus the spectral sequence collapses and we have
$$
H^*(S( \xi_{ \pi}); \mathbf{F}_p) \approx H^*(B \pi; \mathbf{F}_p)
\oplus \tilde{H}^*(S( \xi^{ \times p}); \mathbf{F}_p)^{ \pi}.
$$
Also if $q \neq p$ we have $H^*(S ( \xi_{ \pi}); \mathbf{F}_q)
\approx H^*(S( \xi^{ \times p}); \mathbf{F}_q)^{ \pi}.$

Let $X_{ \pi }$ be a stable summand of $BU(1)^p$ defined as follows.
The action of $ \pi$ on $BU(1)^p $ induces an action of $ \pi$ on the stable
decomposition of $BU(1)^p$ as a wedge of all smash products
of length $1,\ldots,p-1$, say $Y_{ \pi}$, and a smash product of length $p$.
Then choose $X_{ \pi}$ such that $ NX_{ \pi} = Y_{ \pi} $, where
$N=1+t+\cdots +t^{p-1}$. By the stable equivalence
\begin{equation}
S( \xi^{ \times p}) \rightarrow BU(1)^p \rightarrow Y_{ \pi} \label{eq:161}
\end{equation}
we can consider $X_{ \pi}$ as a stable summand of $S( \xi^{ \times p})$.
For any choice of $X_{ \pi}$, consider the composition of stable maps
\begin{equation}
f_{ \pi}: S( \xi_{ \pi}) \rightarrow X^p_{h\pi} \stackrel{Tr}{ \rightarrow}
E \pi \times BU(1)^p \rightarrow BU(1)^p \rightarrow X_{ \pi}.
\end{equation}
We have to show that the stable map $ \varphi_0 \vee f_{ \pi} $
induces an isomorphism in cohomology for any group of coefficients $\mathbf{F}_q$,
$q$ a prime, and hence gives a stable equivalence by the stable
Whitehead lemma. It follows from the above arguments that
$$
\tilde{H^*}(S( \xi_{ \pi}); \mathbf{F}_p)=
\varphi_0^* \tilde{H^*}(B \pi; \mathbf{F}_p) \oplus Tr_0^*
\tilde {H^*}(S( \xi^{ \times p}); \mathbf{F}_p),
$$
and $ \tilde{H^*}(S( \xi_{ \pi}); \mathbf{F}_q)=
Tr_0^* \tilde{ H^*}(S( \xi^{ \times p}); \mathbf{F}_q),$
when $q \neq p$. The restriction of $ Tr_0$ on $X_{ \pi}$ induces a monomorphism
on $Im Tr_0^*$ since by the transfer property (iv), $ \rho_0^* Tr_0^* = N$
and the restriction of $N$ on $ \tilde{H^*}( X_{ \pi}; \mathbf{F}_q) $ is a monomorphism.
Hence $ (\varphi_0 \vee Tr_0|X_{ \pi })^*$ is an isomorphism and so is
$ (\varphi_0 \vee f_{ \pi})^* $ by the commutative diagram
\begin{equation}
\label{eq:171}
\begin{diagram}
 S( \xi_{ \pi}) &\rTo &  X^p_{h\pi}     \\
\dTo{Tr_0} &  &  \dTo{Tr} \\
 S( \xi^{ \times p})    &   \rTo &  E \pi \times BU(1)^p
\end{diagram}
\end{equation}
This proves Lemma \ref{LemH}.\endproof

\proof[Proof of Theorem \ref{ThmB}]

(a)\qua We consider the restriction of any element $y \in MU^*(X^p_{h\pi})$
to $MU^*(S( \xi_{ \pi}))$. By Lemma  \ref{LemH}
we see
this restriction has the form $\varphi_0^*(u) + f_{\pi}^*(w)$ for some $ u \in \tilde{MU}^*(B \pi) $,
$w \in MU^*(X_{\pi})$. Since the composition
$S( \xi_{ \pi}) \rightarrow X^p_{h\pi} \stackrel{ \varphi }{ \rightarrow} B \pi$
coincides with $ \varphi_0 $, $\varphi^*(u)$ also restricts to $\varphi_0^*(u)$.
By diagram (\ref{eq:171}) there is an element $v \in MU^*(BU(1)^p)$ such that
$Tr^*(v)$ restricts to $f_{\pi}^*(w)$. By exactness
$$
y= \varphi^*(u) + Tr^*(v) + y_1 c_p,
$$
for some  $y_1 \in MU^*(X^p_{h\pi})$. For use in the proof of (c) we observe that (\ref{eq:133}),
(\ref{eq:161}) imply $v$ can be chosen in the direct summand
$MU^*[[x_1,\ldots,x_p]]/(x_1 \cdots x_p)$. Thus
we can assume this expression for $y$ is unique and if $v\neq 0$ then $Tr^*(v)$ restricts
nontrivally in $MU^*(S(\xi_{\pi}))$.

Then suppose $cy=0$. We know that $ \rho^*( \theta) = \mathbf{C} $, hence
$ \rho^*(c)=0 $ and $c Tr^*(v)=Tr^*( \rho^*( c)v)=0$
by Frobenius reciprocity.
So we have $c \varphi^*(u) +  cy_1 c_p=0. $
We want to prove $\varphi^*(u) \in ImTr^*$.
Applying $i^*$ we have $0 = i^*(c\varphi^*(u)) + i^*(cy_1c_p) = zu$ since
$i^*\varphi^* = id$, $c = \varphi^*(z)$, and $i^*(c_p)=0$.
Hence $u \in Ann(z) = ImTr^*_{\mathbf{Z}/p}$.
By naturality of the transfer $\varphi^*(u) \in ImTr^*$.
Thus $c\varphi^*(u)=0$ and therefore $cy_1c_p=0$.
Multiplication by $c_p$ is injective by Lemma  \ref{LemH}, hence $cy_1=0$.
Since $ dim(y_1) = dim(y)- 2p $
iterating this argument gives us statement (a).

(b)\qua This follows from the fact that the right
homomorphism in the short exact sequence from Lemma \ref{LemH} is
multiplication by the Euler class $c_p( \xi_{ \pi})$.

(c)\qua Let $y \in Ker \rho^* $. Since $ \varphi \rho = * $
we have
$$
0= \rho^*(Tr^*(v)) + \rho^*(y_1 c_p).
$$
If the first summand is not zero it restricts nontrivially in $MU^*(S(\xi^{\times p}))$ by
definition of $v$. However the second summand restricts to zero since
$\rho^*(c_p)=x_1 \cdots x_p$,
\ $\rho^*(y_1c_p)= \rho^*(y_1)x_1\cdots x_p$ and $x_1\cdots x_p$ restricts to zero as
the Euler class.  Hence both summands are zero. Furthermore multiplication by
$x_1\cdots x_p$ is a
monomorphism hence $\rho^*(y_1)=0$. So
$$
y= \varphi^*(u)+y_1 c_p= \varphi^*(u)+( \varphi^*(u_1)+y_2c_p)c_p=
\varphi^*(u)+ \varphi^*(u_1)c_p + y_2c_p^2.
$$
Repetition of this process proves (c).

(d)\qua
The fact that $c, c_1, c_2 $ multiplicatively
generate $MU^*B(\pi \wr U(1))$ follows from Lemma \ref{LemH}.
The relations $c_1 = c_1^*$, $c_2 = c_2^*$ follow from the bundle relation
$$
\xi_{ \pi} \otimes_{\mathbf{C}} \theta = ( \xi \otimes_{\mathbf{C}}
\rho^*( \theta))_{ \pi} = \xi_{ \pi},
$$
which in turn follows from transfer property (i).

So we have to prove that the Chern classes $c, c_1, c_2 $ with these relations are
a complete system of generators and relations.
Let us use the
splitting principle to write formally
$$ \xi_{\pi} = \eta_1 + \eta_2;$$
$$ u_1 = c_1( \eta_1); $$
$$ u_2 = c_1( \eta_2).$$
Let $F(x,y)= \sum \alpha_{ij} x^i y^j$
be the formal group law.
Using the bundle relation above and applying the Whitney formula for the first and
second Chern classes, we obtain two relations of the form:
$$
F(u_1,c)+F(u_2,c)=c_1
$$
and
\begin{equation}
F(u_1,c)F(u_2,c)=c_2;
\end{equation}
or in terms of $ c, c_1 = u_1+ u_2, c_2 = u_1 u_2 $
\begin{equation}
F(u_1,c)+F(u_2,c)-c_1=c(2+ \sum \beta_{ijk} c^i c_1^j c_2^k)=0 \label{eq:19}
\end{equation}
and
\begin{equation}
F(u_1,c)F(u_2,c)-c_2=c(c_1+ \sum \gamma_{ijk} c^i c_1^j c_2^k)=0, \label{eq:20}
\end{equation}
for some coefficients $ \beta_{ijk} $ , $ \gamma_{ijk} \in MU^*(pt)$.

We claim that relations (\ref{eq:19}) and (\ref{eq:20}) are equivalent to the
following two
obvious transfer relations for $Tr^* : MU^*[[x,tx]] \to MU^*(B( \pi \wr U(1)))$
$$
cTr^*(1)=0
$$
and
$$
cTr^*(x)=0.
$$
Rewrite relations (\ref{eq:19}) and (\ref{eq:20}) as follows:
$$
ca = 0, \ \ \  where \ \ a = 2 + \alpha_{11}c_1 + \sum_{k \geq 2} \alpha_{k1}(u_1^k+u_2^k)
+o(c),
$$
$$
cb = 0, \ \ \ where \ \ b = c_1+2 \alpha_{11}c_2+ \sum_{k \geq 2} \alpha_{k1}(u_1^{k-1}+
u_2^{k-1})c_2+o(c),
$$
and the $\alpha_{ij}$ are the coefficients of the formal group law.

By the first part of Theorem 3.2, $a \in ImTr^*$. Also by transfer property (vi)
$$
\rho^*(a) = \rho^*(Tr^*(1)+ \alpha_{11} Tr^*(x)+
\sum_{k \geq2} \alpha_{k1} Tr^*(x^k)).
$$
\noindent Thus by Lemma 2.1
$$
Tr^*(1)+ \alpha_{11} Tr^*(x)+ \sum_{k \geq2} \alpha_{k1} Tr^*(x^k)=
(F(u_1 ,c)+F(u_2 ,c)-c_1)/c;
$$
similarly $b \in ImTr^*$ and
$$
Tr^*(x)+ \alpha_{11} Tr^*(1)c_2+ \sum_{k \geq2} \alpha_{k1} Tr^*(x^{k-1})c_2=
(F(u_1 ,c)F(u_2 ,c)-c_2)/c.
$$
Now since
$$
x^k = x^{k-1}(x+tx)-x^{k-2}(xtx),
$$
transfer property (i) and the computation of $Tr^*(x)$ is sufficient
for the computation of $Tr^*(x^k)$, $k \geq2$ (see also Corollary 3.6, Remark 3.7).
So we have
\begin{equation}
Tr^*(1)(1+g_0)+Tr^*(x)h_0=(F(u_1 ,c)+F(u_2 ,c)-c_1)/c, \label{eq:21}
\end{equation}
and
\begin{equation}
Tr^*(1)g_1+Tr^*(x)(1+h_1)=(F(u_1 ,c)F(u_2 ,c)-c_2)/c,  \label{eq:22}
\end{equation}
where $g_0, h_0, g_1, h_1 \in MU^*(B(\pi \wr U(1)))$.
This proves (d).

This completes the proof of Theorem \ref{ThmB}.\endproof

Formula (\ref{eq:22}) for computing $Tr^*(x)$ is complicated;
let us give a simpler form. Consider again (\ref{eq:20}).
Note that the coefficient $ \gamma_{00k} \in MU^*(pt)$ contains a factor 2:
the element
$$
c_1+ \sum \gamma_{ijk} c^i c_1^j c_2^k
$$
annihilates $c$ and hence belongs to $ImTr^*$. On the other hand
$$
\rho^* Tr^* =1+t; \:  \rho^*(c)=0; \:  \rho^*(c_1) = x+tx; \: \rho^*(c_2)= xtx,
$$
hence applying $ \rho^*$ we have that
$$
x+tx + \sum \gamma_{0jk} (x+tx)^j (xtx)^k
$$
belongs to $Im(1+t)$. So $ \gamma_{00k}(xtx)^k=2 \gamma_k (xtx)^k$,
that is, $ \gamma_{00k}= 2 \gamma_k$ for some coefficient $ \gamma_k$.

Recall that on the other hand $ F(c,c)=0 $ that is $2c=o(c^2)$. So
$ \gamma_{00k} c =o(c^2)$, hence taking into account the relation $F(c,c)=0$
we can rewrite (\ref{eq:20}) after division by
$$
1+ \sum_{i,k \geq0} \gamma_{i1k} c^i c_2^k
$$
(the coefficient at $cc_1 )$ as follows
\begin{equation}
cc_1=d_0 c + d_2 cc_1^2 +\cdots +d_n cc_1^n +\cdots , \label{eq:23}
\end{equation}
where $d_k=d_k(c,c_2) \in MU^*[[c,c_2]]$ and $d_0(0,c_2)=0$;
the lower index $n$ indicates the coefficient at $cc_1^n$.
Since
$$ \rho^*(Tr^*(x)-c_1) =0,
$$
it follows from Theorem \ref{ThmB}(c) that there exist elements
$$
\delta_j \in \tilde{MU}^*(B\pi)
$$
such that
$$
Tr^*(x)= c_1+ \sum_{j \geq0} \varphi^*( \delta_j) c_2^j.
$$
Using the inclusion $i: B\pi \rightarrow B(\pi \wr U(1))$ we have
$$
i^*(c_1)=i_0^*(c); \: i^*Tr^*(x)=0; \: i^*(c_2)=0,
$$
thus
$$
\varphi^*( \delta_0 ) =-c.
$$
For the calculation of the other elements $ \delta_j $ recall that $cTr^*(x)=0$,
hence
\begin{equation}
cc_1^n =-c \varphi^*(\delta^n); n \geq1, \label{eq:24}
\end{equation}
where
$$
\delta=-c+ \sum_{j \geq1} \delta_j c_2^j.
$$
Combining (\ref{eq:23}) and (\ref{eq:24}), we have the following:
\begin{Prop}
\label{PropXX} The elements $\delta_j $, $j>0$ can be determined from the recurrence relations
which arise from the following formula in $MU^*(B\pi)[[c_2]]$
$$
\delta= d_0 + \sum_{i \geq 2} d_i \delta^i.
$$
\end{Prop}
\proof
By definition the element $ \delta- d_0 - \sum_{i \geq 2} d_i \delta^i $
belongs to $Ker \rho^*$. On the other hand this element is
annihilated by $c$ hence
$$
\delta- d_0 - \sum_{i \geq 2} d_i \delta^i \in ImTr^* \bigcap Ker( \rho^*) =0
$$
by Lemma 2.1.
\endproof

\begin{Cor}
\label{CorHH} For the elements
$ \delta_j \in \tilde{MU}^*(B\mathbf{Z}/2)$, constructed
in \ref{PropXX}, the following formula holds
in $MU^*B(\mathbf{Z}/2 \wr U(1))$

$$
Tr^*(x)= c_1- c+ \sum_{j \geq1} \varphi^*( \delta_j ) c_2^j .
$$
\end{Cor}

In fact we have proved Theorem \ref{ThmJ} for $p=2$. The general case, analogous but
more technical, is given next.

\subsection{Proof of Theorem \ref{ThmJ}}

Note that by the definition of $\omega_k$ the difference
$Tr^*(\omega_k)-c_k$ is an element of $Ker(\rho^*)$. Thus Theorem 3.2 (c) implies
existence of the elements $ \delta_i^{(k)} $ in Theorem \ref{ThmJ}.

First let us elucidate the meaning of the relations
$$
\xi_{ \pi} \otimes_{\mathbf{C}} \theta = \xi_{ \pi}
$$
in the general case of $B(\pi \wr U(1))$.

Again, we can use the splitting principle and write formally
$$
\xi_{ \pi} = \eta_1 + \eta_2 +\cdots + \eta_p ; \  u_m = c_1( \eta_m), \ m=1,\ldots,p.
$$
Applying the Whitney formula for the relation
$$
\eta_1 \otimes_{\mathbf{C}} \theta +\cdots + \eta_p \otimes_{\mathbf{C}} \theta = \eta_1 +\cdots + \eta_p,
$$
and taking into account that $c_m = c_m( \xi_{ \pi})$ is the elementary
symmetric function $ \sigma_m (u_1,\ldots,u_p)$  we have
\begin{equation}
\sigma_m(F(u_1 , c),\ldots, F(u_p , c))= c_m, \label{eq:25}
\end{equation}
$m=1,\ldots,p$, or in terms of $ c , c_1,\ldots,c_p$ we have
$$
 c (p+ \sum \beta^0_{i_0,i_1,\ldots,i_p} c^{i_0}c_1^{i_1}\cdots c_p^{i_p})=0;
$$
and
\begin{equation}
  c((p-k)c_k+
\sum \beta^k_{i_0,i_1,\ldots,i_p} c^{i_0}c_1^{i_1}\cdots c_p^{i_p})=0;
\label{eq:26}
\end{equation}
for $k=1,\ldots,p-1$ and some
$ \beta^0_{i_0,i_1,\ldots,i_p}$,$ \: \beta^k_{i_0,i_1,\ldots,i_p} \in MU^*(pt)$.

We claim that these relations are
equivalent to the following obvious relations
$$
 cTr^*(1)=0,
$$
and
$$
  cTr^*(\omega_k)=0,
$$
for the elements $\omega_k \in MU^*(BU(1))^p$, $k=1,\ldots,p-1$
defined above.

For the proof of our claim multiply the $k$-th relation from (\ref{eq:26})
by $p_k= (p-k)^{-1}$ in $\mathbf{F}_p$.
Then by Theorem \ref{ThmB}, $Ann(c)$ coincides with $ImTr^*$  hence
(\ref{eq:25}) implies that
$$
p_k( \sigma_{k+1}(F(u_1 , c),\ldots, F(u_p , c))-
c_{k+1})/c =Tr^*(a_k),
$$
for some $a_k$ which we have to find.  Let us write
$$
\rho^*(p_k( \sigma_{k+1}(F(u_1 , c),\ldots, F(u_p , c))-
c_{k+1})/ c)=
g^{(k)}( \sigma_1,\ldots, \sigma_p )
$$
$$
= \sigma_k (1+ g_k^{(k)}( \sigma_1,\ldots, \sigma_p))+
\sum_{j \neq k, 1 \leq j \leq {p-1} }  \sigma_j g_j^{(k)}
( \sigma_j, \sigma_{j+1},\ldots, \check { \sigma_k },\ldots,\sigma_p)
$$
$$
=N(\omega_k)(1+ g_k^{(k)}( \sigma_1,\ldots, \sigma_p))+
\sum_{j \neq k, 1 \leq j \leq {p-1} } N(\omega_j)g_j^{(k)}
( \sigma_j, \sigma_{j+1},\ldots, \check { \sigma_k },\ldots,\sigma_p).
$$
Here the symbol $ \check { \sigma_k }$ indicates absence of the corresponding
term. So we have
$$
p_k (\sigma_{k+1}(F(u_1 , c),\ldots, F(u_p , c))-c_{k+1})/c
$$
$$
=Tr^*(\omega_k)(1+ g_k^{(k)}(c_1,\ldots,c_p))+
\sum_{j \neq k, 1 \leq j \leq {p-1} } Tr^*(\omega_j) g_j^{(k)}
(c_j ,c_{j+1},\ldots, \check {c_k },\ldots,c_p) ,
$$
and
$$
[ \sigma_1(F(u_1 , c),\ldots, F(u_p , c))-c_{1}]/c
$$
$$
=Tr^*(1)(1+g_0^{(0)}(c_1,\ldots,c_p))+
\sum_{ 1 \leq j \leq {p-1} } Tr^*(\omega_j) g_j^{(0)}(c_j ,c_{j+1},\ldots,c_p).
$$
This proves our claim.

For computing $\delta_i^{(k)}$ we start with the equations (\ref{eq:26}) and
rewrite them as
\begin{equation}
  cf_k(c ,c_1,\ldots,c_p)=0, \quad k=1,\ldots,p-1. \label{eq:27}
\end{equation}
These are equations in a power series
algebra $MU^{*}(B\pi)[[c_p]]$, since we know
$cc_k \in cMU^{*}(B\pi)[[c_p]].$

We now want to find explicitly formal series
\begin{equation}
 \delta^{(k)}(c_p)=\sum_{i\geq 0}\delta _i^{(k)}(c ){c_p}^i \label{eq:28}
\end{equation}
such that
$$ Tr^{*}(\omega_k)=c_k+\delta ^{(k)}(c_p) $$
and hence
\begin{equation}
 c c_k^j=- c(\delta ^{(k)}(c_p))^j,\quad j\geq 1. \label{eq:29}
\end{equation}
For this we want to replace  the equations (\ref{eq:27}) by the equations
\begin{equation}
 c \widetilde{f}_k(c ,\delta ^{(1)}(c_p),\ldots,\delta^{(p-1)}(c_p),c_p)=0, \label{eq:30}
\end{equation}
where $\widetilde{f}_k\in Ker \rho_{ \pi} ^*$ is a series whose coefficient at $\delta^{(k)}$
is invertible.
In fact $\widetilde{f_k}=0$ since we know that
$Ann(c )=ImTr^{*}$ and $ Ker( \rho^*)\cap ImTr^{*}=0$ by Lemma 2.1.

Then  equating each coefficient of the resulting series
\begin{equation}
{g}_k(c_p) =
\widetilde{f}_k(c ,\delta ^{(1)}(c_p),\ldots,\delta^{(p-1)}(c_p),c_p)=0
\label{eq:31}
\end{equation}
in the ring $MU^{*}(B\pi)[[c_p]]$ to zero we will obtain $p-1$ infinite strings of equations in
$MU^{*}(B\pi)$. Assuming $\delta _i^{(l)}$ are already found for $i<n$ we get
\begin{equation}
\delta _n^{(k)}=\psi_{n,k}((\delta _i^{(1)})_{i \leq n}, \ldots,(\delta
_i^{(k)})_{i< n},\ldots, (\delta_i^{(p-1)})_{i\leq n}), \label{eq:32}
\end{equation}
a system of linear equations in $\delta _n^{(l)}$, $l=1,\ldots,p-1$
with invertible determinant and coefficients in $MU^*[[c]]$.
Since the $\delta_0^{(l)}$ are already known as
$l$-th Chern classes of the bundle $1+\theta +\cdots + \theta ^{p-1}$,
by induction on $n$ we can solve formally (25) to get
\begin{equation}
\delta _n^{(k)}(c)=\tilde{\psi}_{n,k}((\delta_i^{(l)})_{i<n}). \label{eq:321}
\end{equation}
This gives $\delta_n^{(k)}= \delta_n^{(k)}(z)\in MU^*(B\pi)$ obviously
satisfying our equations.

Now for the remaining equation (23) we proceed as follows: let us look at the
term $cf_k(0,0,\ldots,0,c_p)$ in equations (20). Note that $f_k(0,0,\ldots,0,c_p)$
is divisible by $p$:

$f_k\in Ann(c )=ImTr^{*}\Rightarrow \
\rho_{\pi} ^{*}f_k\in ImN\Rightarrow
f_k(0,\ldots,0,\sigma _p)\in ImN \Rightarrow f_k(0,\ldots,0,\sigma_p)$
is divisible by $p$.

Next using the relation $[p]_F(c )=0$ we know that $pc$ is
divisible by $ c^2;$ hence each occurrence of $pc $ in these
equations can be replaced by terms with higher powers of $c$. So
$ c f_k(0,0,\ldots,0,c_p)$ can be replaced by a term divisible by
$c^2$.

Also the $k$-th relation from (20) contains the term
$c(p-k)c_k$,
and for the condition (24) we have to multiply the $k$-th equation
from (20) by $(p-k)^{-1}$, the inverse of $p-k$ in $\mathbf{F}_p$,
and as above we can replace $ c(p-k)c_k$ by $ cc_k$+(terms divisible
by $ c^2)$. Then we use (22) and substitute the series
$\delta ^{(k)}$ in the resulting equations, thus obtaining (23).

This completes the proof for $E=MU$, which is the universal example of complex
oriented cohomology theories. From this result we can descend to all $E$. 
\endproof

We now turn to computation of $Tr^*$ in general.

\begin{Cor}
\label{CorX} For all primes $p$,
Theorem \ref{ThmJ} enables us to explicitly compute the transfer
homomorphism for those
polynomials $a \in \tilde{MU}^*[[x_1,\ldots,x_p]]$ for which
$Na=a+ta+\cdots +t^{p-1}a$ is symmetric in $ x_1,\ldots,x_p $.
\end{Cor}
\proof
If
$Na= \sigma_1 a_1( \sigma_1,\ldots, \sigma_p)+ \ldots+
\sigma_{p-1}a_{p-1}( \sigma_1,\ldots, \sigma_p),$ then
$$
Tr^*(a)=Tr^*(\omega_1)
a_1(c_1,\ldots,c_p) +\cdots + Tr^*(\omega_{p-1})a_{p-1}(c_1,\ldots,c_p)).
$$
To see this let $ \hat a = \omega_1 a_1( \sigma_1,\ldots,
\sigma_p) +\cdots + \omega_{p-1} a_{p-1}( \sigma_1,\ldots, \sigma_p)$. Then  $N(a-
\hat a)=0$, that is, $a- \hat a \in Im(1-t)$,  hence
$Tr^*(a)=Tr^*( \hat a)$.
\endproof

\noindent{\bf Remark 3.7}\qua For $p=2$ one has recurrence formulas for $Tr^*(x^k)$, $k \geq1$.
$$
Tr^*(x) = Tr^*(\omega_1)
$$
$$
Tr^*(x^k) = Tr^*(x^{k-1})c_1 -Tr^*(x^{k-2})c_2
$$
This follows using the formula $x^k = x^{k-1}(x+tx)-x^{k-2}(xtx)$.

\section {Transferred Chern classes for $\Sigma_p$-coverings}

If we consider a $p$-local complex oriented cohomology $E^*$ then by standard
transfer arguments (see Lemma 4.3 below) $E^*(B\Sigma_p)$ is isomorphic to the
subring of $E^*(B\pi)$
invariant under the action of the normalizer of $\pi$ in $\Sigma_p$. The results of this
section imply the elements $\delta^{(k)}_{i} \in \tilde{E}^*(B\pi)$ from
Theorem 3.1 are invariant under this action. This defines elements
$\tilde{ \delta}^{(k)}_i \in \tilde{E}^*(B\Sigma_p)$ which we use for computing the transfer.

In this section we consider $BP^*(X^p_{h \Sigma_p})$ for
$X = \mathbf{C}P^{\infty}$
and for the covering projection
$$
{\rho}_{\Sigma_p} : E{\Sigma_p} \times X^p \to X^{p}_{h\Sigma_p}
$$
we give a formula for the transfer homomorphism
\begin{equation}
{Tr_{\Sigma_p}}^*: BP^*(X^p) \to BP^*(X^p_{h\Sigma_p})
\label{eq: e}
\end{equation}
using the elements $\tilde{\delta}^{(k)}_i$.

\subsection{Brown-Peterson cohomology
of $(\mathbf{C}P^{\infty})^p_{h \Sigma_p}$.}

We need definitions analogous to those of Section 2, with the cyclic
group replaced by the symmetric group. The $p$-fold product, ${\xi}^{\times p}$,
of the canonical line bundle over $X^p$ extends to an
$p$-dimensional bundle
\begin{equation}
\label{eq: bbb}
\xi_{\Sigma_p} =E\Sigma_p \times_{ \Sigma_p} \xi^{ \times p}
\end{equation}
over $X^p_{h\Sigma_p}$ classified by the inclusion
$X^p_{h\Sigma_p} = B(\Sigma_p \wr U(1))\hookrightarrow BU(p)$.
Let $c_i = c_i(\xi_{ \Sigma_p})$.
Then ${\rho_{\Sigma_p}}^*(c_i) = c_i({\xi}^{\times p}) =\sigma_i$,
the $i$-th symmetric polynomial in the $x_j$, where
$BP^*(X^p)=BP^*[[x_1,\ldots,x_p]].$

Then we have the projection
\begin{equation}
\varphi : X^p_{h\Sigma_p} \rightarrow B\Sigma_p,
\end{equation}
induced by the factorization $\Sigma_p \wr U(1)/U(1)^p=\Sigma_p$
and the inclusion
\begin{equation}
i : B\Sigma_p \rightarrow X^p_{h\Sigma_p},
\end{equation}
induced by the inclusion of $\Sigma_p$ in $\Sigma_p \wr U(1)$.

\begin{Def}
\label{DefB}
Let $\tilde{c}_i = {Tr_{\Sigma_p}}^*(x_1x_2 \cdots x_i)$ for
 $i = 1, \ldots, p-1.$
\end{Def}

\begin{Lem}
\label{LemE}
${\rho_{\Sigma_p}}^*(\tilde{c}_i)=i!(p-i)!\sigma_i$.
\end{Lem}

\proof
${\rho_{\Sigma_p}}^*(\tilde{c}_i) = {\rho_{\Sigma_p}}^*Tr_{\Sigma_p}^*(x_1x_2\cdots x_i)=
N_{\Sigma_p}(x_1x_2 \cdots x_i)$.
For each
subset of $i$ integers
$\{j_1, j_2,\ldots, j_i \}$ with $1 \leq j_k \leq p$,
there are $i!$ bijections $\{1, 2, \ldots , i \} \to
\{j_1, j_2,\ldots, j_i\}$ and $(p-i)!$ bijections $\{i+1, i+2, \ldots, p\} \to
\{1,2,\ldots,p\}\backslash \{j_1, j_2,\ldots, j_i\}$. Thus there are $i!(p-i)!$ summands
of $x_{j_1}x_{j_2}\cdots x_{j_i}$ in $N_{\Sigma_p}(x_1x_2\cdots x_i)$. \endproof

We  recall
$$
    BP^*(B\pi) = BP^*[[z]]/([p]z)
$$
with $|z| = 2$. The corresponding computation for $BP^*(B\Sigma_p)$ is also
known \cite{S}. For the reader's convenience we derive the result in a form
useful for our purposes.
\begin{Lem}
\label{LemHS}
As a $BP^*$ algebra
$$
{\rm(i)}  \hskip .5in  BP^*(B\Sigma_p) = BP^*[[y]]/(yTr^*_{ \Sigma_p}(1)),
$$
\noindent where $y$ and $Tr^*_{ \Sigma_p}(1)$ are  uniquely determined by
$ \rho_{ \pi, \Sigma_p}^*(y) = z^{p-1}$ and
$$
{\rm(ii)}  \hskip .7in    \rho^*_{ \pi, \Sigma_p}(Tr^*_{ \Sigma_p}(1))=
(p-1)!Tr^*_{ \pi}(1)=(p-1)![p](z)/z.
$$
\noindent In particular  $|y| = 2(p-1)$.
\end{Lem}

\proof
(ii)\qua Applying the double coset formula (transfer property (vi)) to
$$
 BP^*(Be)  \overset{Tr_{e, \Sigma_p}^*}\longrightarrow BP^*(B\Sigma_p)
 \overset{\rho_{\pi, \Sigma_p}^*}\longrightarrow BP^*(B\pi),
$$
the statement follows from Quillen's formula (1).

\noindent (i)\qua The relation $yTr^*_{\Sigma_p}(1) = 0$ is a consequence
of Frobenius reciprocity.
To see that it is the defining relation we recall that
the cohomology of $B\Sigma_p$
with simple coefficients in $\mathbf{Z}_{(p)}$ is
$$
    H^*(B\Sigma_p; \mathbf{Z}_{(p)}) = \mathbf{Z}_{(p)}[y]/(py)
$$
where $|y| = 2(p-1)$. This follows easily from the mod-$p$ cohomology
and the Bockstein spectral sequence.

Also $H^*(B\pi;  \mathbf{Z}_{(p)}) = \mathbf{Z}_{(p)}[z]/(pz)$ where
$|z| = 2$. The map $\rho_{\pi, \Sigma_p} : B\pi \to B\Sigma_p$ yields
$\rho_{\pi, \Sigma_p}^*(y) = x^{p-1}$.

Now the Atiyah-Hirzebruch-Serre spectral sequence for
$BP^*(B\Sigma_p)$ is
$$
E_2 = H^*(B\Sigma_p; BP^*) =
BP^*[y]/(py)  \Longrightarrow BP^*(B\Sigma_p).
$$
Since $y$ is even dimensional, the sequence collapses at
$E_2 = E_{\infty}$. Thus\break $BP^*(B\Sigma_p)$ is generated by
$y$ as a $BP^*$ algebra.

For the group $W = N_{\Sigma_p}(\pi)/\pi \approx \mathbf{Z}/(p-1)$,
$|W|$ is prime to $p$, hence by the standard transfer argument
$\rho_{\pi, \Sigma_p}^*: BP^*(B\Sigma_p) \to BP^*(B\pi)$ is an
injective map of $BP^*$ algebras.
Since $\rho^*_{ \pi, \Sigma_p}(yTr^*_{ \Sigma_p}(1)) = p!z^{p-1}$ plus
terms of higher filtration, $yTr^*_{\Sigma_p}(1) = 0$ is the only relation.
\endproof

Relating $\pi$ and $\Sigma_p$ we have a lift of
$\rho_{\pi, \Sigma_p}$
\[
\begin{diagram}
X^{p}_{h\pi} &\rTo{\tilde{\rho}_{\pi , \Sigma_p}} & X^{p}_{\Sigma_p}  \\
\dTo{\varphi} & & \dTo{\varphi} \\
B\pi & \rTo{{\rho}_{\pi, \Sigma_p}} & B\Sigma_p
\end{diagram}
\]

\begin{Lem}
\label{LemEE}
$\tilde{\rho}_{\pi,\Sigma_p}^{*}(\widetilde{c}_k) = k!(p-k)!
Tr_{\pi}^{*}(\omega_k).$
\end{Lem}

\proof
Note that modulo $Im(1-t)$ we have $\Sigma
g^{*}(x_1x_2\cdots x_k)=k!(p-k)!\omega_k$ summed over $\Sigma_p /\pi $.
Applying the double coset formula
\begin{align*}
\rho_{\pi, \Sigma_p}^{*}(\widetilde{c}_k) =
\rho_{\pi, \Sigma_p}^{*}Tr_{\Sigma _p}^{*}(x_1x_{2}\cdots x_k)&=\\
Tr_\pi^{*}\sum_{g\in \Sigma_p/\pi}& g^{*}(x_1x_2\cdots x_k)=
k!(p-k)!Tr_\pi ^{*}(\omega_k).\tag*{\qed}
\end{align*}

Let $c=\varphi^*(y) \in BP^{2(p-1)}(X^p_{h\Sigma_p}).$

\begin{Lem} $ImTr_{\Sigma_p}^*$ is contained in the $BP^*$
algebra generated by $$c,\tilde{c}_1,\ldots ,\tilde{c}_{p-1},c_p.$$
\end{Lem}

\proof By the K\"unneth isomorphism,
\begin{equation}
  BP^*(X^p) =  BP^*(X)^{\otimes p} = F \oplus T
\label{eq:c}
\end{equation}
as a $\pi$-module, where $F$ is free and  $T$ is trivial. Explicitly
a $BP^*$ basis for $T$ is $\{x_{1}^i \cdots x_{p}^i,i\geq 0\}$, while
a $BP^*$ basis for $F$ is $\{x_{1}^{i_1} \cdots x_{p}^{i_p},
i_j \geq 0\}$
where not all the exponents are equal.

By Lemma \ref{LemHS} $Tr_{\Sigma_p}^*(1)$ is
a power series in $c$.
Now recall from \cite{Hb}, p. 44,
that we can consider $BP^*(X^p)$,  $X=\mathbf{C}P^{\infty}$ as a
free $BP^*[[\sigma_1,\ldots ,\sigma_p]]$ module generated by $1$ and the elements
$x_{1}^{i_1} \cdots x_{p}^{i_p} \in F$, with
$ 0 \leq i_j \leq p-j $.
So by Frobenius reciprocity it suffices to compute the transfer on these
monomials.
Summed over the symmetric group
$ \sum g^*(x_{1}^{i_1} \cdots x_{p}^{i_p})$ is a symmetric function and hence has the form
$$
\sum_{ \pi} \sum_{ \Sigma_p / \pi} g^*(x_{1}^{i_1} \cdots x_{p}^{i_p}) =
\sigma_1 s_1 +\cdots + \sigma_{p-1}s_{p-1}=
\sum_{ \pi} (\omega_1 s_1+\cdots + \omega_{p-1}s_{p-1}),
$$
for the elements $\omega_k$ from Theorem 3.1  and symmetric functions
$s_1,\ldots ,s_{p-1}$.
Hence modulo $ker N_{\pi} =Im(1-t)$, $t \in \pi$, we have the following
equation in
$F$
$$
\sum_{ \Sigma_p / \pi} g^*(x_{1}^{i_1}  \cdots x_{p}^{i_p}) =
\omega_1 s_1 +\cdots + \omega_{p-1}s_{p-1}.
$$
The left sum consists of $(p-1)! $ elements each having
the same transfer value. Also $\omega_k $ is the sum of $p^{-1} \binom pk$
elements $x_{i_1}\cdots x_{i_k}$; on each of these elements the transfer evaluates
to $Tr_{ \Sigma_p}^* (x_1 \cdots x_k) = \tilde{c}_k$.
Thus Frobenius reciprocity and Lemma 4.2 is all
that is needed for computing $Tr^*_{\Sigma_p}$.
\endproof

Recall the elements $\delta^{(k)}_i \in \tilde{BP}^*(B \pi)$ derived from Theorem \ref{ThmJ}
by naturality.
By the standard transfer argument again the map induced by
$\tilde{ \rho}_{\pi, \Sigma_p} : X^p_{h\pi} \rightarrow X^p_{h\Sigma_p} $, the
lift of $ \rho_{\pi, \Sigma_p} : B\pi \rightarrow B\Sigma_p$,
is also injective. Moreover for $BP^*(X^p_{h\Sigma_p})$ the ring structure
is completely determined by the following:

\begin{Thm}
\label{ThmBPalg}
As a $BP^*$ algebra
$$
BP^*(X^p_{ h \Sigma_p})=
BP^*[[c, \tilde{c}_1,\ldots, \tilde{c}_{p-1}, c_p]]
/(cTr^*_{ \Sigma_p}(1), c \tilde{c}_i )
$$
and one has the formula
$$
\widetilde{c}_k-k!(p-k)!c_k=\Sigma _{i\geq 0}\varphi^{*}(\tilde{\delta}
^{(k)}_i)c_p^i,\quad k=1,\ldots,p-1, \
$$
where the elements  $ \tilde{ \delta}^{(k)}_i \in \tilde{BP}^*(B \Sigma_p)$
are determined by
$$
\rho_{ \pi, \Sigma_p}^{*}( \tilde{ \delta}^{(k)}_j)=k!(p-k)! \delta^{(k)}_j,
\quad   j\geq 0.
$$
\end{Thm}

For the proof we follow that of Theorem 3.2. Let
$$
S( \xi_{ \Sigma_p} )=E \Sigma_p \times_{ \Sigma_p} S(\xi^{ \times p})
$$
be the sphere bundle of the bundle $ \xi_{\Sigma_p}$ of (28).
$X^p_{h \Sigma_p}$ is homotopy equivalent to the disk bundle
$D( \xi_{ \Sigma_p})=E \Sigma_p \times_{ \Sigma_p} D( \xi^{ \times p})$.
Then we have the obvious inclusion
$i_0 : B \Sigma_p \rightarrow S( \xi_{ \Sigma_p})$
and projection
$ \varphi_0 : S( \xi_{ \Sigma_p}) \rightarrow B \Sigma_p$
with fiber $ S( \xi^{ \times p})$. $ \varphi_0 i_0 $ is the identity. Thus
stably $ B\Sigma_p$ is a wedge
summand of $S( \xi_{ \Sigma_p})$. As for the other summand let
$$
X_{ \Sigma_p}= \vee_{i=1}^{p-1} E\Sigma_i {\times}_{\Sigma_i} BU(1)^{ \wedge i}.
$$
By the standard transfer argument, localized at $p$, $X_{\Sigma_p}$ is a stable summand
of $\vee_{i=1}^{p-1} E\Sigma_i \times BU(1)^{ \wedge i}$ and hence of
$E\Sigma_p \times BU(1)^{\times p}$. From this we derive the following result.

\begin{Lem} One has a stable equivalence localized at $p$
$$
\varphi_0 \vee f_{ \Sigma_p} : S( \xi_{ \Sigma_p}) \to B \Sigma_p
\vee X_{ \Sigma_p},
$$
with $f_{ \Sigma_p} $, the composition of stable maps
$$
f_{ \Sigma_p}: S( \xi_{ \Sigma_p }) \rightarrow X^p_{h \Sigma_p}
\stackrel{Tr_{ \Sigma_p }}{ \rightarrow}
E \Sigma_p \times BU(1)^p \rightarrow X_{ \Sigma_p }.
$$
\end{Lem}

\proof
The inclusion $i_0$ splits off $ \varphi^*_0 H^*(B\Sigma_p)$ in
$H^*(S( \xi_{ \Sigma_p}))$. Furthermore in mod-$p$ cohomology
$$
H^*(S( \xi^{ \times p}))=
\mathbf{F}_p[x_1,\ldots,x_p]/ (\sigma_p),
$$
hence
$$
 H^*(S( \xi^{ \times p}))^{\Sigma_p}
 \approx \mathbf{F}_p[ \tilde{c}_1 ,\ldots, \tilde{c}_{p-1}],
$$
by Lemma 4.2.

Then $ H:=\tilde{H}^*(S( \xi^{ \times p}))$ is a free $\pi$ module and
$H^*(\Sigma_p; H) \subseteq H^*(\pi; H).$ Thus
\begin{displaymath}
    H^*(\Sigma_p; H) = H^{\Sigma_p}\ \ \  if \ \ \ *=0
\end{displaymath}
\begin{displaymath}
     \ \ \ \ \ \ \ \ \ \ \ \ \ \ \    = 0 \ \ \ \ \ \   if \ \ \ \ *>0.
\end{displaymath}
Therefore there is an isomorphism
$$
H^*(S( \xi_{ \Sigma_p})) \overset{\rho^* \oplus i_0^*}
\longrightarrow \tilde{H}^*(S( \xi^{ \times p}))^{\Sigma_p}
\oplus H^*(B\Sigma_p)
$$
where $\rho: S( \xi^{ \times p}) \to S( \xi_{ \Sigma_p})$ is the projection.
We have to prove that the first summand
is $f^*_{ \Sigma_p} \tilde{H}^*(X_{ \Sigma_p})$.

By naturality of the transfer we have the commutative diagram
\[
\begin{diagram}
S( \xi^{ \times p})      & \rTo{} & BU(1)^p    & \rTo{} X_{ \Sigma_p} \\
\uTo{Tr_{0}}     &       & \uTo{Tr_{\Sigma_p}}  \\
S( \xi_{ \Sigma_p})      & \rTo{} & X^p_{h \Sigma_p}
\end{diagram}
\]
Thus $f_{ \Sigma_p}$ coincides
with $ \tilde{f}_{ \Sigma_p}$, the map $Tr_0$ followed by the horizontal maps
in the above diagram. We wish to show the restriction of $Tr_0^*$ to the image of
$H^*(X_{\Sigma_p})$ is an isomorphism onto $H^*(S( \xi^{ \times p}))^{\Sigma_p}$.

Now considering the transfer for the $\Sigma_i$ coverings
$$
   E\Sigma_i \times BU(1)^{\wedge i} \to E\Sigma_i \times_{\Sigma_i} BU(1)^{\wedge i},
$$
it follows from transfer properties (ii) and (vi) that
$H^*(E\Sigma_i \times_{\Sigma_i} BU(1)^{\wedge i})$ is a submodule of
$H^*(E\Sigma_i \times BU(1)^{\wedge i})$ generated by $\Sigma_i$ norms
of monomials in $x_1, x_2,\ldots, x_i$, with non-increasing degrees. From this it is
straightforward that $H^*(X_{\Sigma_p})$ and $H^*(S( \xi^{ \times p}))^{\Sigma_p}$
have the same ranks in each dimension. Thus we are reduced to showing the desired
map is injective.

However,
for any monomial $x$ in $x_1, x_2,\ldots, x_i$, we have
$$
    Tr_0^*(N_{\Sigma_i}(x)) = i!Tr_0^*(x)
$$
by naturality of the transfer.
Thus the restriction of $Tr_0^*$ to the image of $H^*(X_{\Sigma_p})$ will be
a monomorphism if $Tr_0^*$ is non-zero on polynomials consisting of monomials
with non-increasing degrees. This in turn will follow if the norm $N_{\Sigma_p}$
is non-zero on such polynomials. In fact we claim: 1) $N_{\Sigma_p}$ is non-zero
on any monomial $x^I = x_{1}^{i_1}\cdots  x_{p-1}^{i_{p-1}}$,
and 2) different monomials with
non-increasing degrees in $x_1, \ldots ,x_i$, $i<p$ are in different $\Sigma_p$ orbits.

Claim 2) is clear. To see 1) let $J = (j_1,\ldots,j_p)$ and
$x^J = x_{1}^{j_1} \cdots  x_{p}^{j_p}$, all of whose exponents are not equal.
Then we will show the coefficient of $x^J$ in
$N_{\Sigma_p}(x^J)$ is prime to $p$.
The isotropy subgroup of $x^J$ is the finite product
$\Sigma_{n_1} \times \Sigma_{n_2} \times \cdots  < \Sigma_p$
where $n_j$ is the number of terms of $J$ equaling $j$.
This group has order ${n_1}!{n_2}!\cdots $ which is prime to $p$.
Hence $N_{\Sigma_p}(x) = ({n_1}!{n_2}!\cdots )x^J$ + other monomials proving
the claim.

Thus $ \varphi_0 \vee f_{ \Sigma_p}$ induces an isomorphism and
hence is a $p$-local stable equivalence.
\endproof

This implies the following:

\begin{Lem} The long exact sequence for the pair
$(D( \xi_{ \Sigma_p}), S( \xi_{ \Sigma_p}))$ gives the
following short exact sequence
$$
0 \leftarrow BP^* (S( \xi_{ \Sigma_p})) \leftarrow BP^*(X^p_{h \Sigma_p}) \leftarrow
BP^*((X^p_{h \Sigma_p})^{ \xi_{ \Sigma_p}}) \leftarrow 0.
$$
\end{Lem}

Indeed the left arrow is an epimorphism by Lemma 4.7 and hence the right
arrow is a monomorphism.
\endproof

Now the proof of Theorem \ref{ThmBPalg} is completely analogous to that of
Theorem \ref{ThmB} taking into account additionally that any element
$y \in BP^*(X_{h\pi}^p)$ has the form
$$
y= {\varphi}^*(u) +g( \tilde{c}_1,\ldots, \tilde{c}_{p-1})+y_1 c_p
$$
for some $u \in BP^*(X_{h\pi}^p)$, where $g$ denotes some formal power series and
$y_1 \in BP^*(X_{h\pi}^p)$. This follows by  Lemma 4.5 and Lemma 4.8.
\endproof

\section{Calculation of the elements $\delta^{(k)}_i$ and
$\tilde{\delta}^{(k)}_i$ in Morava $K$-theory }

In this section we work in Morava $K$-theory $K(s)^*$
and give an alternative, better algorithm  for explicit computations.

Fix a prime $p$ and an integer $s \geq 0$, then $K(s)^* = \mathbf{F}_p[v_s, {v_s}^{-1}]$
with $|v_s| = -2(p^s -1)$.
By a result of W\"urgler \cite{Wg} there is no restriction on $p$:
although $K(s)$ is not a commutative ring spectrum for $p=2$,
we shall consider only those spaces whose Morava $K$-theory is even dimensional. This
implies the deviation from commutativity is zero.

We recall
$$
   K(s)^*(B\pi) = K(s)^*[z]/(z^{p^s})
$$
where $|z| = 2$.

As in Lemma 4.3. we have:

\begin{Lem}
\label{LemD}
{\rm(i)}\qua $\rho_{\pi, \Sigma_p} : B\pi \to B\Sigma_p$ induces an isomorphism
of $K(s)^*$ algebras
$$
\rho_{\pi, \Sigma_p}^*: K(s)^*(B\Sigma_p) \overset{\approx}\longrightarrow \{K(s)^*(B\pi)\}^W,
$$
where $W = N_{\Sigma_p}(\pi)/\pi \approx \mathbf{Z}/(p-1)$.
Computing invariants yields
$$
  K(s)^*(B\Sigma_p) = K(s)^*[y]/(y^{m_s}),
$$
where $\rho_{\pi, \Sigma_p}^*(y) = z^{p-1}$ and $m_s = [(p^s - 1)/(p-1)]+1$.

{\rm(ii)}\qua  $Tr_{\Sigma _p}^{*}(1)=-v_s{y}^{m_s-1}.$

\end{Lem}

Then combining Theorem 4.6 and Remark 2.2 we have

$$
K(s)^*(X^p_{h\Sigma_p}) = K(s)^*[[c, \tilde{c}_1,\ldots,\tilde{c}_{p-1}, c_p]]/
(c^{m_s}, c\tilde{c}_i).
$$
Our main result in this section is the following:

\begin{Prop}
\label{PropE} We can construct explicit elements
$\delta^{(k)}_i \in K(s)^*(B\pi)$ such that

{\rm(1)}\qua In $K(s)^*(X^p_{h \pi})$ the following formula holds
$$
c_k( \xi_{ \pi})=Tr_{ \pi}^*(\omega_k)-
\sum_{0 \leq i \leq p^s} \varphi_{ \pi}^*( \delta^{(k)}_i)c_p^i
( \xi_{ \pi}).$$

{\rm(2)}\qua In $K(s)^*( X^p_{h \Sigma_p})$ one has
$$
c_k( \xi_{ \Sigma_p})=Tr_{\Sigma_p}^*(x_1\cdots x_k)-
\sum_{0 \leq i \leq p^s} \varphi_{ \Sigma_p}^*( \tilde{\delta}^{(k)}_i)
c_p^i( \xi_{ \Sigma_p }),
$$
with  $\tilde{\rho}_{\pi,\Sigma_p}^{*}( \tilde{ \delta}^{(k)}_i) =
 k!(p-k)! \delta^{(k)}_ik.$

{\rm(3)}\qua The value of $Tr_{ \pi}^*(c_1( \xi_i))$ is determined by
$$
c_1( \xi_{ \pi})=Tr_{ \pi}^*(c_1( \xi_i))+
v_s \sum_{1 \leq j \leq s-1}c^{p^s-p^j}c_p^{p^{j-1}}( \xi_{ \pi})
$$
where $\xi_i$ is the pullback of the canonical line bundle $\xi$ by projection
$BU(1)^p \to BU(1)$ on the $i$-th factor.
\end{Prop}

We are grateful to D.~Ravenel for supplying us with the proof of the
following result.

\begin{Lem}
\label{LemF1}
For the formal group law in Morava K-theory $K(s)$, $s>1$, we have

$$
F(x,y) \equiv x + y - v_s \sum_{0<j<p}p^{-1}\binom{p}{j}(x^{p^{s-1}})^{j}(y^{p^{s-1}})^{p-j}
$$
modulo $x^{p^{2(s-1)}}$ (or modulo $y^{p^{2(s-1)}}$).
\end{Lem}

\proof
This result can be derived from the recursive formula for the FGL
given in 4.3.9 \cite{R}.  For the FGL in Morava $K$-theory
it reads
$$
F(x,y) = {\sum_{i\geq 0}}^F v_s^{e_i} w_i(x,y)^{p^{i(s-1)}}
$$
where $w_i$ is a certain homogeneous polynomial of
degree $p^i$ defined by 4.3.5 \cite{R} and $e_{i}= (p^{is}-1)/ (p^{s}-1)$.
In particular $\omega_{0}=x+y$,
\begin{displaymath}
w_{1}=-\sum_{0<j<p}p^{-1}\binom{p}{j}x^{j}y^{p-j} ,
\end{displaymath}
\noindent and $w_{i} \notin (x^{p},y^{p})$.

We find it more convenient to express $F(x,y)$ as
$$
F(x,y) = F(x+y, v_s w_1(x,y)^{p^{s-1}}, v_s^{e_{2}}
w_2(x,y)^{p^{2 (s-1)}},\ldots).
$$
Then for $s>1$ we can reduce modulo the ideal $v_{s}^{e_{2}}(x^{p^{2
(s-1)}}, y^{p^{2 (s-1)}})$ and get
\begin{align*}
F(x,y) &\equiv F ( x + y, v_s w_1(x,y)^{p^{s-1}})\\
    & =
F (x + y+ v_s w_1(x,y)^{p^{s-1}},
   v_{s}w_{1}(x + y, v_s w_1(x,y)^{p^{s-1}})^{p^{s-1}}, \dots )\\
 &\equiv  F (x + y+ v_{s}w_1(x, y)^{p^{s-1}},
        v_{s}w_{1}(x^{p^{s-1}} +
y^{p^{s-1}}, v_{s}^{p^{s-1}} w_1(x,y)^{p^{2(s-1)}}) ),
\end{align*}
and modulo $v_{s}^{1+p^{s-1}}(x^{p^{2 (s-1)}}, y^{p^{2
(s-1)}})$ we have
$$
F(x,y) \equiv x + y + v_s w_1(x,y)^{p^{s-1}}. \eqno{\qed}
$$
Let us write for short $\sigma_k=\sigma_k(x,F(x,z),\ldots,F(x,(p-1)z))$.

\begin{Cor}
The following formula holds in  $K(s)^*(BU(1) \times B\pi)$
$$
\sigma_k=
- \sum_{0 \leq i \leq p^s} \lambda^{(k)}_i
\sigma_p^i+p^{-1} \binom{p}{k} x^kv_sz^{p^s -1},
$$
where $ \lambda^{(k)}_i= \lambda^{(k)}_i(z^{p-1})$ are polynomials in
$z^{p-1}$ and $ \lambda^{(j)}_0=0$, $j=1,\ldots,p-2$,
$\lambda^{(p-1)}_0=-z^{p-1}.$
\end{Cor}

\noindent{ \it{Proof.}} For $1 \leq k \leq p-1$, equating the coefficients of $x^{ip}$, $1 \leq i \leq p^s$ gives a system of linear equations with invertible
matrix of the form $Id + nilpotent$.
Thus the elements $ \lambda^{(k)}_1,\ldots,\lambda^{(k)}_{p^s}$ can be defined as the solution of this system. Of course equating the coefficients at $x^i$ for
$i\neq p,2p,\ldots,p^{s+1}$ will produce other equations in $\lambda^{(k)}_j$, $j=1,\ldots,p^s$.
But these equations
are derived from the old equations above. These additional
equations make the matrix upper triangular. \endproof

Now, let us prove Proposition 5.2 and show that one necessarily
has $\delta^{(k)}_i=\lambda^{(k)}_i$, $i=0,\ldots,p^s$ for
$\lambda^{(k)}_i$ encountered in Corollary 5.4. Thus by Lemma 5.1,
$\delta^{(k)}_i$ is invariant under the action of $W$ and we can
define $\tilde{\delta}^{(k)}_i$ by
$\tilde{\rho}_{\pi,\Sigma_p}^{*} ( \tilde{ \delta}^{(k)}_i) =
k!(p-k)! \delta^{(k)}_i.$

The diagonal map $ \Delta : BU(1) \rightarrow BU(1)^p$ induces an inclusion
$B\pi \times BU(1) \rightarrow X^p_{h \pi}$ and the commutative diagram
\[
\begin{diagram}
E\pi \times BU(1)  & \rTo{1 \times \Delta}  & E\pi \times BU(1)^p  \\
\dTo{ \pi \times 1} &  &  \dTo{ \rho_{ \pi }}  \\
B\pi \times BU(1) & \rTo{                }  & X^{p}_{h\pi}
\end{diagram}
\]
Then $ (1 \times \Delta )^* (\omega_k)= p^{-1} \binom{p}{k} x^k$, $x=c_1( \xi)$.
Hence by transfer properties (i) and (iv) we have for the transfer
$Tr=Tr( \pi \times 1)$:
$$
Tr^*((1 \times \Delta)^* (\omega_k))=p^{-1} \binom{p}{k} x^k Tr^*(1)=
p^{-1} \binom{p}{k} x^k v_s z^{p^s -1}.
$$
On the other hand by the existence of the elements $ \delta^{(k)}_i$
(Theorem 3.1) we have
$$
Tr^*((1 \times \Delta)^*(\omega_k))=
$$
$$
\sigma_k(x,F(x,z),\ldots,F(x,(p-1)z))
+ \sum_{i \geq 0}\delta^{(k)}_i\sigma_p^i(x,F(x,z),\ldots,F(x,(p-1)z)):
$$
$\xi_{\pi}$ restricts to $\sum_i \xi \otimes \theta^i$ on $BU(1) \times B\pi$, thus $c_k(\xi_{\pi})$
to $\sigma_k(x,F(x,z),\ldots,$ $F(x,(p-1)z))$; by Lemma 5.3 and the fact that $z^{p^s}=0$, $[i]z$
may be replaced by $iz$. By Corollary 5.4
$$
\sigma_k(x,F(x,z),\ldots,F(x,(p-1)z))=
$$
$$
- \sum_{0 \leq i \leq p^s } \lambda^{(k)}_{i}\sigma_p^i(x,F(x,z),\ldots,F(x,(p-1)z))
+p^{-1} \binom{p}{k} x ^kv_sz^{p^s -1}.
$$
Then the restriction of $ (1 \times \Delta )^*$ to $Ker \rho^*$ is a
monomorphism \cite{Hun}. This proves
Proposition 5.2.1) and shows
$\delta^{(k)}_i = \lambda^{(k)}_i$ for $0\leq i \leq p^s$ and zero
otherwise. Statement 2) follows from Lemma 4.4. Then 3) follows from the
following explicit formula for $ \sigma_1$:

\begin{Lem} In  $K(s)^*(B \pi \times BU(1))$ one has
$$
\sigma_1=v_s\left(z^{p^s-1}x
+\sum_{i=1}^{s-1}z^{p^s-p^i}\sigma_p^{p^{i-1}}\right).
$$
\end{Lem}

\proof One has
\begin{align*}
\sigma_1&=x+F(x,z)+\cdots +F(x,(p-1)z)=x+x+z+v_s w_1(x^{p^{s-1}},z^{p^{s-1}})\\
&\qquad\qquad\qquad+\cdots +x+(p-1)z+v_s w_1(x^{p^{s-1}},((p-1)z)^{p^{s-1}})\\
&=px+\frac{p(p-1)}2z+v_s\left(\sum_{i=1}^{p-1} w_1(x,iz)\right)^{p^{s-1}}\\
&=v_s\left(\sum_{i=1}^{p-1}
\sum_{j=1}^{p-1}-p^{-1}\binom pji^jx^{p-j}z^j\right)^{p^{s-1}}\\
&=v_s\left(\sum_{j=1}^{p-1}-\left(\sum_{i=1}^{p-1}i^j
\right)p^{-1}\binom pjz^jx^{p-j}\right)^{p^{s-1}}.
\end{align*}
Now
$ \sum_{i=1}^{p-1}i^j$ is an integral linear combination of
$\sigma_k(1,2,\ldots,p-1)$ with $k\le i$, hence by it
is zero for
$i<p-1$ and for $i=p-1$ it is $p-1$.

Thus
\begin{equation}
\sigma_1=-v_s\left((p-1)p^{-1}\binom p{p-1}z^{p-1}x^{p-(p-1)}\right)^{p^{s-1}}=
v_sz^{p^s-p^{s-1}}x^{p^{s-1}}.  \label{eq:11}
\end{equation}
ow since $F(x,z)^p=x^p+z^p$, one has $\sigma_p^p=
\left(x(x+z)\cdots (x+(p-1)z)\right)^p$.

\noindent But again we have $x(x+z)\cdots (x+(p-1)z)=x^p-xz^{p-1}$. Substituting this one
obtains
\begin{align*}
v_s(z^{p^s-1}x&+\sum_{i=1}^{s-1}z^{p^s-p^i}\sigma_p^{p^{i-1}})\\
&=v_s\left(z^{p^s-1}x+z^{p^s-p}\sigma_p
+\sum_{i=2}^{s-1}z^{p^s-p^i}(x^p-z^{p-1}x)^{p^{i-1}}\right)\\
&=v_s\left(z^{p^s-1}x+z^{p^s-p}\sigma_p+\sum_{i=2}^{s-1}z^{p^s-p^i}
(x^{p^i}-z^{(p-1)p^{i-1}}x^{p^{i-1}})\right).
\end{align*}
But it is straightforward to see that
$$\sum_{i=2}^{s-1}z^{p^s-p^i}(x^{p^i}-z^{(p-1)p^{i-1}}x^{p^{i-1}})
=z^{p^s-p^{s-1}}x^{p^{s-1}}-z^{p^s-p}x^p.$$
Hence one has
\begin{align}
v_s (z^{p^s-1}x&+\sum_{i=1}^{s-1}z^{p^s-p^i}\sigma_p^{p^{i-1}})\notag\\
&=v_s\left(z^{p^s-1}x+z^{p^s-p}\sigma_p+z^{p^s-p^{s-1}}x^{p^{s-1}}-z^{p^s-p}x^p
\right). \label{eq:12}
\end{align}
Now one has

$
z^{p^s-p}F(x,kz)=z^{p^s-p}(x+kz+v_s w_1(x^{p^{s-1}},(kz)^{p^{s-1}}))=
z^{p^s-p}(x+kz),
$

\noindent hence
$
z^{p^s-p}\sigma_p=
z^{p^s-p}x(x+z)\cdots (x+(p-1)z)=z^{p^s-p}(x^p-z^{p-1}x)
$

Substituting this into (\ref{eq:12})
gives
$$
v_s\left(z^{p^s-1}x+\sum_{i=1}^{s-1}z^{p^s-p^i}\sigma_p^{p^{i-1}}\right)=
v_sz^{p^s-p^{s-1}}x^{p^{s-1}},
$$
which is $\sigma_1$ by (\ref{eq:11}).
\endproof

We now compute some of the elements
$ \delta^{(k)}_i$ and $ \tilde{\delta}^{(k)}_i$.

First recall from  \cite{H}, \cite{R} that generators for
$$
\begin{array}{ccc}
\pi _{*}BP & \subset  & H_{*}BP \\
\updownarrow  &  & \updownarrow  \\
\mathbf{Z}_{(p)}[v_1,v_2,\ldots] & \subset  & \mathbf{Z}_{(p)}[m_1,m_2,\ldots]
\end{array}
$$
$$
\left| v_n\right| =2(p^n-1)=\left| m_n\right|
$$
are given by
$$v_n= pm_{n} - \sum_{i=1}^{n-1}m_iv_{n-i}^{p^i}.$$
Given a formal group law over a graded ring $R_{*},$
$$F(x,y)= \sum_{i,j} \alpha _{ij}^Rx^iy^j\in R_{*}[[x,y]],\qquad \alpha _{ij}^R\in
R_{2(i+j-1)}$$
there is a ring map $g:MU_{*}\longrightarrow R_{*}$ which induces the formal
group law; that is $g^{*}(\alpha _{ij}^{MU})=\alpha _{ij}^R.$

We use also the following well known formulas
$$
F(x,y)=exp(log x+log y) \ \ \ \ \  and  \ \ \ \ \ log x = \sum_{n \geq0} m_n x^{n+1}
$$
for computing the elements $ \delta_i$ in $BP$ theory by the algorithm
of Section 3.

{\bf Example 1}\qua For $ \delta_1 \in
BP^*(B\mathbf{Z}/2)=BP^*[[z]]/([2](z))$ we have
modulo $z^{8}$:

$\delta_1=
v_1^2 z^2 + (v_1^3 + v_2)z^3 + v_1 z^4 + (v_1^6 + v_1^3 v_2)z^6 + (v_1^4 v_2
+ v_2^2 + v_3) z^7$.
\

Next we give some
results of calculations in Morava $K$-theories, where the formulas are more
tractable. In the following examples
$ \delta^{(k)}_i$ coincides with the coefficient at $ \sigma_p^i$ in the
expression for $\sigma_k$ from Corollary 5.4, $y=z^{p-1}$, and
$ \tilde{\delta} ^{(k)}_i=k!(p-k)! \delta^{(k)}_i$.

\
\

\noindent{\bf Example 2}\qua $p=3, s=2$

$\sigma_{{1}}=v_{{2}}{y}^{3}\sigma_{{3}}+v_{{2}}{y}^{4}x$.

$\sigma_{{2}}=2\,{v_{{2}}}^{2}{y}^{3}{\sigma_{{3}}}^{4}+2\,v_{{2}}{y}^{2}{\sigma_{{3}}}^{2}
+v_{2}{x}^{2}{y}^{4}+2\,y$.

\
\

\noindent{\bf Example 3}\qua $p=5, s=3$

$\sigma_{{1}}=v_{{3}}{y}^{25}{\sigma_{{5}}}^{5}+v_{{3}}{y}^{30}\sigma_{{5}}+v_{{3}}{y}^{31}x$.

$\sigma_{{2}}=4\,{v_{{3}}}^{2}{y}^{25}{\sigma_{{5}}}^{30}+4\,{v_{{3}}}^{2}{y}^{30}{\sigma_{{5}}}^{26}+3\,v_{{3}}{y}^{19}{\sigma_{{5}}}^{10}+v_{{3}}{y}^{24}{\sigma_{{5}}}^{6}+
3\,v_{{3}}{y}^{29}{\sigma_{{5}}}^{2}+2\,v_{{3}}{y}^{31}{x}^{2}$.

$\sigma_{{3}}=2\,{v_{{3}}}^{3}{y}^{25}{\sigma_{{5}}}^{55}+
2\,{v_{{3}}}^{3}{y}^{30}{\sigma_{{5}}}^{51}+{v_{{3}}}^{2}{y}^{19}{\sigma_{{5}}}^{35}+
2\,{v_{{3}}}^{2}{y}^{24}{\sigma_{{5}}}^{31}+{v_{{3}}}^{2}{y}^{29}{\sigma_{{5}}}^{27}+
2\,v_{{3}}{y}^{13}{\sigma_{{5}}}^{15}+v_{{3}}{y}^{18}{\sigma_{{5}}}^{11}+
v_{{3}}{y}^{23}{\sigma_{{5}}}^{7}+2\,v_{{3}}{y}^{28}{\sigma_{{5}}}^{3}+2\,v_{3}{x}^{3}{y}^{31}$.

$\sigma_{{4}}=4\,{v_{{3}}}^{4}{y}^{25}{\sigma_{{5}}}^{80}+4\,{v_{{3}}}^{4}{y}^{30}{\sigma_{{5}}}^{76}+4\,{v_{{3}}}^{3}{y}^{19}{\sigma_{{5}}}^{60}+3\,{v_{{3}}}^{3}{y}^{24}{\sigma_{{5}}}^{56}+
4\,{v_{{3}}}^{3}{y}^{29}{\sigma_{{5}}}^{52}
+4\,{v_{{3}}}^{2}{y}^{13}{\sigma_{{5}}}^{40}+2\,{v_{{3}}}^{2}{y}^{18}{\sigma_{{5}}}^{36}+
2\,{v_{{3}}}^{2}{y}^{23}{\sigma_{{5}}}^{32}+4\,{v_{{3}}}^{2}{y}^{28}{\sigma_{{5}}}^{28}+
4\,v_3{y}^{7}{\sigma_{{5}}}^{20}+v_3{y}^{12}{\sigma_{{5}}}^{16}+
4\,v_3{y}^{17}{\sigma_{{5}}}^{12}+v_3{y}^{22}{\sigma_{{5}}}^{8}+
4\,v_3{y}^{27}{\sigma_{{5}}}^{4}+v_3{x}^{4}{y}^{31}+4\,y$. \endproof

\section{ Transfer and  $K(s)^*(X^{p}_{h\Sigma_p})$}

Let $X$ be a CW complex whose Morava $K$-theory
$K(s)^*(X)$  is even dimensional and finitely generated as a module over
$K(s)^*$.

In this section we study the transfer homomorphism in this more general context. We extend some
results of Hopkins-Kuhn-Ravenel \cite{HKR} to spaces. We
consider the Atiyah-Hirzebruch-Serre (later abbreviated AHS) spectral sequence:
\begin{equation}
{E_2}^{*,*}(\pi, X) = H^*(\pi; K(s)^*X^p) \Rightarrow K(s)^*(X^{p}_{h\pi}).
\label{eq:a}
\end{equation}
By the K\"unneth isomorphism
\begin{equation}
K(s)^*X^p \overset{\approx}\longrightarrow (K(s)^*X)^{\otimes p} \label{eq:b}.
\end{equation}
Then $K(s)^*X^p$ is a $\pi$ module where $\pi$ acts
by permuting factors (see \cite{HKR}, Theorem\ 7.3).

An element $x \in K(s)^*(X)$ is called {\it good} if there is a finite
cover $Y \to X$ together with an Euler class $y \in K(s)^*(Y)$ such that
$x = Tr^*(y)$ where $Tr^*:K(s)^*(Y) \to K(s)^*(X)$ is the transfer.
The space $X$ is called {\it good} if $K(s)^*(X)$ is spanned over $K(s)^*$
by good elements.

Let $\gamma = \varphi^*(z)$, where
$$
\varphi: X^p_{h\pi} \to B\pi
$$ is the projection
and let  $\{x_j, j \in \mathcal{J} \}$ be a $K(s)^* $ basis for $K(s)^*(X)$.
Hunton \cite{Hun1} has shown that if $K(s)^*(X)$ is concentrated in even dimensions
then so is $K(s)^*(X^p_{h\pi})$. We adopt the stronger hypothesis that $X$ is good
and derive a stronger result,
following the argument of \cite{HKR} Theorem\ 7.3 for classifying spaces.

\begin{Prop}
\label{PropA} Let $X$ be a good space.

{\rm(i)}\qua As a $K(s)^*$ module
$K(s)^*(X^p_{h\pi})$ is free with basis
$$
\{ \ {\gamma}^i \otimes (x_j)^{\otimes p}\ \ | \ \ 0 \leq i < p^s, j \in \mathcal{J}   \}
$$
and
$$
\{\sum_{(i_1, i_2,\ldots, i_p)= I} 1 \otimes x_{i_1}\otimes x_{i_2} \otimes \cdots 
\otimes x_{i_p}\ \ | \ \ I\in \mathcal{P}_p \ \ \}
$$
where  $I = \{(i_1, i_2,\ldots, i_p) \}$ runs over the set $\mathcal{P}_p$
of $\pi$-equivalence classes of
$p$-tuples of indices $i_j \in \mathcal{J}$ at least two of which are not equal.

{\rm(ii)}\qua $X^p_{h\pi}$ is good.
\end{Prop}

\proof (i)\qua
By the K\"unneth isomorphism,

\begin{equation}
   K(s)^*(X)^{\otimes p} = F \oplus T
\label{eq:ck}
\end{equation}
as a $\pi$-module, where $F$ is free and  $T$ is trivial. Explicitly
a $K(s)^*$ basis for $T$ is $\{(x_i)^{\otimes p},i \in \mathcal{J} \}$, while
a $K(s)^*$ basis for $F$ is $\{x_{i_1}\otimes x_{i_2} \otimes \cdots  \otimes x_{i_p},
i_j \in \mathcal{J} \}$
where not all the factors are equal. Then
$$
    H^*(\pi; F) = F^{\pi} \ \ \ \ \  if \ \ \ *=0
$$
$$    \ \ \ \ \ \ \ \ \ \ \ \ \   = 0     \ \ \ \ \ \ \    if \ \ \ *>0
$$
and
$$
    H^*(\pi; T) = H^*(B\pi) \otimes T.
$$
Thus  ${E_2}^{0,*}(\pi, X) = K(s)^*(X^{p}_{h\pi})^{\pi} = F^{\pi} \oplus T$.

To continue the proof we recall the covering projection
$$
   {\rho}_{\pi}: E\pi \times X^p \to X^{p}_{h\pi},
$$
its associated transfer homorphism
\begin{equation}
\label{eq: aaa}
  Tr^* = Tr_{\pi}^*: K(s)^*(X^p) \to K(s)^*(X^{p}_{h\pi}),
\end{equation}
and induced homomorphism
$$
    {\rho_{\pi}}^*: K(s)^*(X^{p}_{h\pi}) \to K(s)^*(X^p).
$$
Similar maps are defined
for the group $\Sigma_p$. Then ${\rho_{\pi}}^* Tr^* = N$, where $N = N_{\pi}$
is the trace map.

Thus we have established the following lemma.

\begin{Lem}
\label{LemA}
If $y \in  K(s)^*(X^p)$ is good then
there exists a good element $z \in K(s)^*(X^{p}_{h\pi})$ such
that $\rho_{\pi}^*(z) = N(y)$. \endproof
\end{Lem}

\begin{Lem}
\label{LemB}
 If $x \in K(s)^*(X)$ is good then there is a good
element $z \in K(s)^*(X^{p}_{h\pi})$ such that $\rho_{\pi}^*(z) = x^{\otimes p}$.
\end{Lem}

\proof By assumption there is a finite covering $f:Y \to X$
and an Euler class $e \in K(s)^*(Y)$ such that $x = Tr^*(e)$.
Now consider the covering
$$
    \phi = f \times \cdots  \times f : Y^p \to X^p
$$
which extends to a covering
$$
    1\times \phi:Y^{p}_{h\pi} \to X^{p}_{h\pi}
$$
and yields a map of coverings
\[
\begin{diagram}
E\pi \times Y^p  & \rTo{{\rho}_{\pi}} & Y^{p}_{h\pi}  \\
\dTo{1 \times \phi } &  &  \dTo{1 \times \phi }  \\
E\pi \times X^p & \rTo{{\rho}_{\pi}} & X^{p}_{h\pi}
\end{diagram}
\]
The class $e^{\otimes p}$ is an Euler class for $Y^p$. Since the transfer
is natural and commutes with tensor products we have
$$
   {\rho_{\pi}^*}Tr^*(1 \otimes e^{\otimes p})= Tr^*{\rho_{\pi}}^*(1 \otimes e^{\otimes p}) =
Tr^*(e^{\otimes p}) = Tr^*(e) \otimes \cdots  \otimes Tr^*(e) = x^{\otimes p}.
\eqno{\qed}$$

\begin{Cor}
\label{CorA} ${E_2}^{0,*}(\pi, X)$ consists of permanent
cycles which are good. \qed
\end{Cor}

Thus as differential graded $K(s)^*$ modules, there is an
isomorphism of spectral sequences
$$
   ({E_r}^{*,*}(\pi,pt) \otimes_{K(s)^*} T) \oplus F^{\pi} \overset{\approx}\longrightarrow
{E_r}^{*,*}(\pi,X).
$$
Thus it follows that as a $K(s)^*$ algebra, $K(s)^*(X^{p}_{h\pi})$
is generated
by\nl $K(s)^*(B\pi)$, $T$, and $F^{\pi}$.

(ii)\qua The proof of \cite{HKR} Theorem\ 7.3 carries over. This completes the
proof of Proposition \ref{PropA}.
\endproof

\noindent{\bf Remarks}\qua (1)\qua From the periodicity of the cohomology
of a cyclic group \cite{CE}
Proposition XII, 11.1, we have isomorphisms
$$
   H^t(\pi; K(s)^*(X^p)) \overset{\cdot{z}} \longrightarrow H^{t+2}(\pi; K(s)^*(X^p))
$$
for $t >0$ and
$$
  H^0(\pi; K(s)^*(X^p))/Im(N)  \overset{\cdot{z}} \longrightarrow H^{2}(\pi; K(s)^*(X^p)).
$$
Thus multiplication by $z$ is also injective
on $T$ at the $E_2$ term.

(2)\qua ${\rho_{\pi}}^*Tr^* = N$, thus modulo $ker({\rho_{\pi}}^*)$ we have
\begin{equation}
Tr^*(x_{i_1}\otimes x_{i_2} \otimes \cdots \otimes x_{i_p})= \sum_{\sigma \in \pi}
1 \otimes x_{\sigma(i_1)}\otimes x_{\sigma(i_2)} \otimes \cdots 
\otimes x_{\sigma(i_p)}.
\label{eq d}
\end{equation}
Note that if the $i_j$ in (\ref{eq d}) are equal, the right hand side is zero.
However
$$
Tr^*({x_j}^{\otimes p}) =
 1 \otimes {x_j}^{\otimes p}\cdot Tr^*(1).
$$
We now turn to $K(s)^*(X^p_{h\Sigma_p})$.

Let $c = \varphi^*(y)$ where $\varphi: E\Sigma_p \times_{\Sigma_p}X^p \to
B\Sigma_p$ is the projection.

\begin{Prop}
\label{PropB} Let $X$ be a good space.
As a $K(s)^*$ module
$K(s)^*(X^p_{h\Sigma_p})$ is free with basis
$$
       \{ \ c^i \otimes (x_j)^{\otimes p}\ \ | \ 0\leq i < m_s, j \in \mathcal{J}  \ \}
$$
and
$$
\{\sum_{(i_1, i_2,\ldots, i_p)= I} 1 \otimes x_{i_1}\otimes x_{i_2} \otimes \cdots 
\otimes x_{i_p}\ \ | \ \ I\in \mathcal{E}_p \ \ \}
$$
where  $I = \{(i_1, i_2,\ldots, i_p) \}$ runs over the set $\mathcal{E}_p$
of $\Sigma_p$-equivalence classes of
$p$-tuples of indices $i_j \in \mathcal{J}$ at least two of which are not equal.
\end{Prop}

\proof Since $|W|$ is prime to $p$, the result follows from the AHS spectral sequence,
as in the proof of Proposition \ref{PropA}.\endproof

\section{Applications}

{\bf 7.1}\qua {$\pi\wr({\mathbf Z}/p^n)$}
\medskip

We now turn to  $G_n = \pi  \wr ( {\mathbf Z}/p^n) $ where $\pi = {\mathbf Z}/p$.
Then $BG_n = X_{h\pi}^p$ for $X = B\mathbf{Z}/p^n$.
Consider the AHS spectral sequence for
$$
  B({\mathbf Z}/p^n)^p \to BG_n \overset{\varphi}\to B\pi.
$$
Then
$$
    E_2^{p,q} = H^*(\pi ; K^*(s)(B({\mathbf Z}/p^n)^p)),
$$
where
$$
K^*(s)(B({\mathbf Z}/p^n)^p)  = (K(s)^*[z]/(z^{p^{ns}}))^{\otimes p} = F \oplus T,
$$
where $F$ and $T$ as in (36) above are free (resp. trivial) $\pi$ modules.

Let $\gamma = {\varphi^*(z)}$ where $ K(s)^*(B\pi) = K(s)^*[z]/(z^{p^s})$ as above.

\begin{Prop}
\label{PropF}
As a $K(s)^*$ module $K(s)^*(BG_n)$ is free
with basis

$$
\{ {\gamma}^i \otimes (z^j)^{\otimes p}\ \ 0 \leq i < p^s,  0 \leq j < p^{ns}  \}
$$
and
$$
\{\sum_{(i_1, i_2,\ldots, i_p)= I} 1 \otimes z^{i_1}\otimes z^{i_2} \otimes \cdots 
\otimes z^{i_p}\ \ | \ \ I\in \mathcal{P}_p(n) \ \ \},
$$
where  $I = \{(i_1, i_2,\ldots, i_p) \}$ runs over the set $\mathcal{P}_p(n)$
of $\pi$-equivalence classes of
$p$-tuples of integers $\{ 0 \leq i_j < p^{ns}\}$ at least two of which
are not equal.
\end{Prop}

\proof This spectral sequence computation is exactly analogous
to that of Proposition \ref{PropA}. \endproof

\noindent{\bf Remarks}\qua (i)\qua For $X = {\mathbf C}P^{\infty}$, Proposition \ref{PropF} gives another
derivation of $K(s)^*(X_{h\pi}^p)$.
Since ${{\mathbf C}P^{\infty}}^{\wedge}_p =
[colim_n \ B({\mathbf Z}/p^n)]^{\wedge}_p$,
we have $K(s)^*(X_{h\pi}^p) = lim_n \ K(s)^*(BG_n)$.

\noindent (ii)\qua $G_n$ is good for $K(s)^*$ by \cite{HKR} Theorem 7.3.

By analogy with Section 6 we have the following:

\begin{Lem}
\label{LemG}

\rm{(i)}\qua $Im(Tr^*)\cdot \gamma = 0.$

\rm{(ii)}\qua $Tr^*(1) = v_s {\gamma}^{p^{s}-1}.$

\rm{(iii)}\qua If $y \in T$ then $Tr^*(y) = y \cdot Tr^*(1).$

\end{Lem}

Finally we consider the case $p=2$ where $c = \gamma$. If $n=1$,
$G_n \approx D_8$, the dihedral group of order 8; the rings $K(s)^*(BD_{2^k})$
were determined by Schuster \cite{Sh},\cite{SY}. In general we have a partial
result:

\begin{Prop}
\label{PropG}
Let $p=2$. As a
$K(s)^*(B\pi)$ algebra, $K(s)^*(BG_n)$
is generated by $\tilde{c}_1, c_2$ subject to the
relation $\tilde{c}_1 \cdot c = 0$
and the relations ${\tilde{c}_1}^{2^{ns}} = {c_2}^{2^{ns}} = 0$ modulo
terms divisible by $c$.
\end{Prop}

\proof $\tilde{c}_1 \cdot c = 0$ by Lemma \ref{LemG}.
The other
relations
hold in
the $E_{\infty}$ term of the
spectral sequence. The only possible extensions are those on the fiber,
involving $\tilde{c}_1, {c}_2$. \endproof

Similar results hold for $\Sigma_p \wr {\mathbf Z}/p$ and
$\Sigma_p \wr \Sigma_p$ for $p$ odd.

\medskip
{\bf 7.2}\qua {$p$-groups with cyclic subgroup of index $p$}.
\medskip

In this section we consider the class of $p$-groups with a (necessarily normal)
cyclic subgroup of index $p$. It is known \cite{Br}, Theorem 4.1, Chapter IV, that every $p$-group of
this form is isomorphic to one of the groups:

(a)\qua $\mathbf{Z}/q$ $(q=p^n, n \geq1)$.

(b)\qua $\mathbf{Z}/q \times \mathbf{Z}/p$ $(q=p^n, n \geq1)$.

(c)\qua $\mathbf{Z}/q \rtimes \mathbf{Z}/p$ $(q=p^n, n \geq 2)$, where the canonical
generator of $\mathbf{Z}/p$ acts on $\mathbf{Z}/q$ as multiplication by $1+p^{n-1}$.
This group is called the modular group if $p \geq3 $ and the quasi-dihedral
group if $p=2, n \geq4 $.

For $p=2$ there are three additional families.

(d)\qua Dihedral $2$-groups $D_{2m}=\mathbf{Z}/m \rtimes \mathbf{Z}/2$, $(m \geq2)$, where
the generator of $\mathbf{Z}/2$ acts on $\mathbf{Z}/m$ as multiplication by $-1$. If $m=2^n$,
$D_{2m}$ is a $2$-group. Note that $D_4$ belongs to (b) and $D_8$ belongs
to (c).

(e)\qua Generalized quaternion $2$-groups. Let $\mathbf{H}$ be the algebra of quaternions
$\mathbf{R} \oplus \mathbf{R}i \oplus \mathbf{R}j \oplus \mathbf{R}k$. For $m \geq2$ the
generalized quaternion
group $Q_{4m}$ is defined as the subgroup of the multiplicative group ${\mathbf{H}}^*$ generated
by $x=e^{ \pi i/m}$ and $y=j$. $\mathbf{Z}/2m$ generated by $x$ is normal and has
index $2$. If $m$ is a power of $2$, $Q_{4m}$ is a $2$-group. In the
extension $0 \rightarrow \mathbf{Z}/2m \rightarrow Q_{4m} \rightarrow \mathbf{Z}/2 \rightarrow 0$
the generator of $\mathbf{Z}/2$ acts on $\mathbf{Z}/2m$ as $-1$.  In particular $Q_8$
is the group of quaternions $ \lbrace \pm 1, \pm i, \pm j, \pm k \rbrace $.

(f)\qua Semi-dihedral groups. $\mathbf{Z}/q \rtimes \mathbf{Z}/2$ $(q=p^n, n \geq3 )$, where the
generator of $\mathbf{Z}/2$ acts on $\mathbf{Z}/q$ as multiplication by $-1+2^{n-1}$.

Consider now the task of computing the stable Euler class,
$Tr_G^*(1)$, for the universal $G$-covering $EG \rightarrow BG $.

For the case (a) there is the well known formula of Quillen (\ref{eq: qq})
$$
Tr_{\mathbf{Z}/q}^*(1)=[q]_F(z)/z
$$
in $MU^*(B\mathbf{Z}/q)=MU^*[[z]]/([q]_F (z))$.

For the case (b) the answer follows from transfer property (ii):
$Tr_G=Tr_{\mathbf{Z}/q} \wedge Tr_{\mathbf{Z}/p} $;

In cases (d),(e) and (f) $Tr_G^*$ is the composition of two transfers
$Tr_{\mathbf{Z}/q}^*$ and $Tr_{C,G}^* : MU^*(BC) \rightarrow MU^*(BG)$, where
$C$ is the corresponding cyclic subgroup. So we have to compute
$Tr^*_{C,G}(z^i)$, $i \geq1 $ and we can apply our results
for $B\mathbf{Z}/2 \wr U(1)$, namely Remark 3.7.

Similarly for the case (c), $Tr_G $ is the composition
$Tr_{\mathbf{Z}/q, G} Tr_{\mathbf{Z}/q}$ and we can apply Corollary 3.6.

This task is trivial for wreath products $\mathbf{Z}/p \wr \mathbf{Z}/p^n$ since
$Tr_{\mathbf{Z}/p^n}^*(1)$ is symmetric in $z_1,\ldots,z_p $ in the ring
$$
MU^*((B\mathbf{Z}/p^n)^p)=MU^*[[z_1,\ldots,z_p]]/([p^n](z_1),\ldots,[p^n](z_p))
$$
and hence invariant under the $\mathbf{Z}/p$ action. So in this case we need
only Quillen's formula.

Finally we note that if $G$ is the modular group of case (c), Brunetti \cite{Bi} has completely
computed
the ring $K(s)^*(BG)$. The relations are quite simple but the generators are technically complicated.
In a future paper we plan to use transferred Chern classes to give a more natural presentation.

\medskip
{\bf 7.3\qua Other examples}
\medskip

Consider the semi-direct products
$G=(\mathbf{Z}/p)^n \rtimes \mathbf{Z}/p$ where the generator $ \alpha$ of $\mathbf{Z}/p$
acts on $H_n=\mathbf{Z}/p[T]/(T^n)$ by $1- \alpha =T$, $1\leq n \leq p$.
Then every $\mathbf{Z}/p[\mathbf{Z}/p]$-module is a direct sum of the modules $H_n$.
As shown by Yagita \cite{Y} and Kriz \cite{Kr}, these semi-direct products are good in the
sense of Hopkins-Kuhn-Ravenel.

We recall
$$
K(s)^*(B(\mathbf{Z}/p)^n)=K(s)^* [[z_1,\ldots,z_n]]/(z_i^{p^s}),
$$
where $z_i$ is the Euler class of a faithful complex line bundle $ \theta_i $
on the $i$-th factor. Then
$\mathbf{Z}/p$ acts on $K(s)^*[z_1,\ldots,z_n]/(z_i^{p^s})$ by
$$
\alpha: z_i \rightarrow F_{K(s)}(z_i,z_{i+1}),\ \ \ \  z_{n+1}:=0,
$$
where $ F_{K(s)}$ denotes the formal group law for Morava $K$-theory.

Our aim is to show how to compute the stable Euler classes in terms of
characteristic classes and the formal group law.

The transfer $Tr^* : K(s)^*EG \rightarrow K(s)^*BG$ is the composition of two
transfers
$$
Tr^*_1 : K(s)^*E((\mathbf{Z}/p)^n) \rightarrow K(s)^*B((\mathbf{Z}/p)^n)
$$
and
$$
Tr^*_2 : K(s)^*B((\mathbf{Z}/p)^n) \rightarrow K(s)^*BG.
$$
Recall also that
$$
Tr_1^*(1)= z_1^{p^s-1} \cdots z_n^{p^s-1}.
$$
It is easy to see that in $K(s)^*((B\mathbf{Z}/p)^n)$ we have
$$
e^{p^s-1}( \alpha^{i_1}( \theta_1))\cdots e^{p^s-1}( \alpha^{i_n}( \theta_1))=
z_1^{p^s-1}\cdots z_n^{p^s-1},
$$
where $e$ is the Euler class and $1 \leq i_1 < \cdots < i_n \leq p$.
Then recall the elements $\omega_n$ from Theorem \ref{ThmJ} and let
$\omega_n(l)$ be the sum of the same monomials after raising to the power $l$.
Since $\omega_n(l)$ consist of $p^{-1} \binom pn$
summands and $p^{-1} \binom pn = (-1)^n/n \mod p$,
we have that in $K(s)^*((BZ/p)^n)$
$$
\eta^*_{ \pi}(\omega_n({p^s-1}))= (-1)^n z_1^{p^s-1}\cdots z_n^{p^s-1}/n,
$$
where the map $ \eta_{ \pi} $, defined in Section 2, sends
$ \xi_i =t^{i-1} \xi_1 $ to $ \alpha^{i-1} \theta_1$. Hence
$$
Tr_G^*(1)=Tr_2^*(Tr_1^*(1))=
Tr_2^*((-1)^n n \eta_{ \pi}^*(\omega_n({p^s-1})))=
$$
$$
(-1)^n n Tr_2^*( \eta_{ \pi}^*(\omega_n({p^s-1}))),
$$
and we have to apply Corollary 3.6.

\Addresses
\recd
\end{document}